\documentclass[10pt]{article}
\usepackage{graphicx}
\usepackage{amssymb}
\usepackage{epstopdf}
\usepackage{amsmath}
\usepackage{amsfonts}
\newcommand{\braket}[2]{\langle #1,#2 \rangle}
\newcommand{\la}{\lambda}

\newcommand{\al}{\alpha}

\def\phi{{\varphi}}

\DeclareSymbolFont{AMSb}{U}{msb}{m}{n}
\DeclareMathSymbol{\N}{\mathbin}{AMSb}{"4E}
\DeclareMathSymbol{\Z}{\mathbin}{AMSb}{"5A}
\DeclareMathSymbol{\R}{\mathbin}{AMSb}{"52}
\DeclareMathSymbol{\Q}{\mathbin}{AMSb}{"51}
\DeclareMathSymbol{\I}{\mathbin}{AMSb}{"49}
\DeclareMathSymbol{\C}{\mathbin}{AMSb}{"43}

\begin{document}

\addtolength{\textheight}{0 cm}
\addtolength{\hoffset}{0 cm}
\addtolength{\textwidth}{0 cm}
\addtolength{\voffset}{0 cm}

\setcounter{secnumdepth}{5}
 \newtheorem{proposition}{Proposition}[section]
\newtheorem{theorem}{Theorem}[section]
\newtheorem{lemma}[theorem]{Lemma}
\newtheorem{coro}[theorem]{Corollary}
\newtheorem{remark}[theorem]{Remark}
\newtheorem{ex}[theorem]{Example}
\newtheorem{claim}[theorem]{Claim}
\newtheorem{conj}[theorem]{Conjecture}
\newtheorem{definition}[theorem]{Definition}
\newtheorem{application}{Application}

\newtheorem{corollary}[theorem]{Corollary}
\def\HADX{{\cal H}_{\rm AD}(X)}
\def\HADY{{\cal H}_{\rm AD}(Y)}
\def\HADH{{\cal H}_{\rm AD}(H)}
\def\HTADX{{\cal H}_{\rm TAD}(X)}
\def\HTADY{{\cal H}_{\rm TAD}(Y)}
\def\HTADH{{\cal H}_{\rm TAD}(H)}
\def\LX{{\cal L}(X)}
\def\LY{{\cal L}(Y)}
\def\LH{{\cal L}(H)}
 \def\ASD{{\cal L}_{\rm AD}(X)}
 \def\ASDY{{\cal L}_{\rm AD}(Y)}
\def\ASDH{{\cal L}_{\rm AD}(H)}
 \def\ASDP{{\cal L}^{+}_{\rm AD}(X)}
  \def\ASDYP{{\cal L}^{+}_{\rm AD}(Y)}
   \def\ASDHP{{\cal L}^{+}_{\rm AD}(H)}
    \def\TADX{{\cal L}_{\rm TAD}(X)}
        \def\TADY{{\cal L}_{\rm TAD}(Y)}
            \def\TADH{{\cal L}_{\rm TAD}(H)}
 \def\CX{{\cal C}(X)}
\def\CY{{\cal C}(Y)}
\def\CH{{\cal C}(H)}

\def\PX{{\cal A}(X)}
\def\PY{{\cal A}(Y)}
\def\PH{{\cal A}(H)}
\def\phi{{\varphi}}
\def\AH{A^{2}_{H}}
\def\B {{\cal B}}

\title{Anti-symmetric Hamiltonians (II): Variational resolutions for Navier-Stokes and other nonlinear  evolutions}
\author{ Nassif  Ghoussoub\thanks{Partially supported by a grant
from the Natural Sciences and Engineering Research Council of Canada.  } \quad  and \quad Abbas Moameni\thanks{Research supported by a postdoctoral fellowship at the University of British Columbia.}
\\
\small Department of Mathematics,
\small University of British Columbia, \\
\small Vancouver BC Canada V6T 1Z2 \\
\small {\tt nassif@math.ubc.ca}\\
 \small {\tt moameni@math.ubc.ca}
\\\\
\date{}
}
 \maketitle

\section*{Abstract}
 The nonlinear selfdual variational principle established in a preceeding paper \cite{G3}
-- though good enough to be readily applicable in many stationary nonlinear partial differential equations --  did not however cover the case of nonlinear evolutions such as the Navier-Stokes equations. One of the reasons  is the prohibitive coercivity condition that is not satisfied by the corresponding selfdual functional on the relevant path space. We show here that such a principle still hold for functionals of the form 
  \begin{equation*}
I(u)= \int_0^T \Big [  L  (t,  u(t),\dot { u}(t)+\Lambda  u(t)) +\langle \Lambda  u(t),  u(t) \rangle \Big ] \, dt +\ell (u(0)- u(T), \frac { u(T)+ u(0)}{2})
\end{equation*}
where $L$ (resp., $\ell$) is an anti-selfdual Lagrangian on state space (resp., boundary space), and $\Lambda$
is an appropriate nonlinear  operator on path space.
As a consequence,  we provide a variational formulation and resolution to evolution equations involving  nonlinear operators such as the Navier-Stokes equation (in dimensions $2$ and $3$) with various boundary conditions.  
In dimension $2$, we recover the well known solutions for the corresponding initial-value problem as well as periodic and anti-periodic ones, while in dimension $3$ we get Leray solutions for the initial-value problems, but also solutions satisfying $u(0)=\alpha u(T)$ for any given $\alpha$ in $(-1,1)$.  Our approach is quite general and does apply to many other situations.

  \section{Introduction} This paper is a continuation of \cite{G3} where the first-named author established a general nonlinear selfdual variational principle, that yields a variational formulation and resolution for several nonlinear partial differential equations which are not normally of Euler-Lagrange type. Applications included   nonlinear transport equations, the stationary Navier-Stokes equations, and the generalized Choquard-Pekar Schr\"odinger equations with certain non-local potentials. The principle also applied to the complex Ginsburg-Landau evolution equations, but could not cover Leray's existence results for Navier-Stokes evolutions in low dimensions. The primary objective of this paper is to develop a sharper   selfdual variational principle to be able to deal with this shortcoming, and to encompass a larger class of nonlinear  evolution equations in its scope of applications. 
  
We first recall the basic concept of  {\it selfduality}.  It relates to the following  class of Lagrangians which play a significant role in our proposed  variational formulation. If $X$ is a reflexive Banach space, and  $L:X\times X^{*} \to \R \cup \{+\infty\}$ is a  convex lower semi-continuous function, that is not identically equal to $+ \infty$,  we say that  $L$ is an {\it anti-selfdual Lagrangian (ASD)}  on $X\times X^*$ if
\begin{equation}
\label{seldual}
L^*( p, x) =L(-x, -p ) \quad \hbox{\rm for all $(p,x)\in X^{*}\times X$}, 
\end{equation}
 where  $L^*$ is the Legendre-Fenchel dual (in both variables) of $L$,  defined on $X^{*}\times X$ as:
 \[
 L^*( q,y)= \sup \{ \braket{q}{x} + \braket{p}{y} - L(x,p);\, x \in X, p \in X^{*}  \}.
 \]
 We shall frequently use the following basic properties of an ASD Lagrangian:
 \begin{equation}\label{obs.1}
\hbox{$ L(x, p)+\langle x, p\rangle\geq 0$   for every $(x, p) \in X\times X^{*}$,}
 \end{equation}
  and the fact that
 \begin{equation}\label{obs.2}
\hbox{  $L(x, p)+\langle x, p\rangle =0$ if and only if $(-p, -x)\in \partial L(x,p)$. }
\end{equation}
We therefore define the  {\it derived vector fields of } $L$ at $x \in X $ to be the -possibly empty- sets
\begin{eqnarray}
\bar \partial L(x):= \{ p \in X^*; L(x,-p)- \langle x,p \rangle=0 \}= \{ p \in X^*;  (p, -x)\in \partial L(x,p) \}.
\end{eqnarray}
These {\it anti-selfdual vector fields} are natural extensions of subdifferentials of convex lower semi-continuous functions. Indeed, the most
basic anti-selfdual Lagrangians are of the form $L(x,p)= \varphi (x)+\varphi^*(-p)$ where  $\varphi$ is such a function in $X$, and $\varphi^*$
is its Legendre conjugate on $X^*,$ in which case  $\bar \partial  L(x)= \partial \varphi (x).$   More interesting  examples of
anti-selfdual Lagrangians are of the form $L(x,p)= \varphi (x)+\varphi^*(-\Gamma x-p)$  where $\varphi$  is a convex and lower semi-continuous
function on $X,$ and  $\Gamma: X\rightarrow X^*$ is a skew adjoint operator. The corresponding anti-selfdual vector field is then 
$\bar \partial L(x)=\Gamma x+ \partial \varphi (x)$.   Actually, it turned out that every  { \it maximal monotone operator}  (we refer to \cite{Br} for this well developed theory) is an {\it anti-selfdual vector field}. This fact proved by the first-named author in \cite{G5} means that ASD-Lagrangians can be seen as the {\it potentials} of maximal monotone operators, in the same way as the Dirichlet integral  is the potential of the Laplacian operator (and more generally as any convex lower semi-continuous energy is a potential for its own subdifferential), leading to a variational formulation and resolution of most equations involving maximal monotone operators.

 In this article,   we develop further the approach  
-introduced in \cite {G3}-  to allow for a variational
resolution of non-linear PDE's of the form 
\begin{equation}\label{obs.3}
\Lambda u+ \bar \partial
L(u)=0,
\end{equation}
and nonlinear evolution equations of the form  
\begin{equation}
\hbox{$\dot u (t)+\Lambda u(t)+\bar \partial L(u(t))
=0$ starting at $u(0)=u_{0}$},
\end{equation}
 where $L$ is an anti-selfdual Lagrangian and
$\Lambda:D(\Lambda)\subset X\to X^{*}$ is a non-linear  {\it regular  map}, that is  if 
 \begin{equation}
 \hbox{$\Lambda$ is weak-to-weak continuous and $u\to \langle \Lambda u, u\rangle$ is weakly lower semi-continuous on $D(\Lambda)$}.
 \end{equation}

We note that 
positive   linear operators are necessarily
{\it regular  maps}, but that there is also 
a wide class of nonlinear regular operators,  such as those appearing in the basic equations of  hydrodynamics and magnetohydrodynamics (see below and \cite{Te}). 

Our approach is based on the following simple observation: If $L$ is an anti-selfdual Lagrangian on $X\times X^*$,
then for any map $\Lambda:D(\Lambda)\subset X\to X^{*}$,   we have from (\ref{obs.1}) and  (\ref{obs.2}) above that 
  \begin{equation}
\label{positive}
I(x):=L(x,\Lambda x)+\langle x, \Lambda x \rangle\geq 0 \quad \hbox{\rm for all $x\in D(\Lambda)$}, 
\end{equation}
and that equation (\ref{obs.3}) is satisfied by $\bar x\in X$ provided the infimum of $I$ is equal to zero and that it is attained at $\bar x$. The following theorem established in \cite{G3} provides conditions under which such an existence result holds. 

\begin{theorem} \label{main.1}  Let $L$ be an anti-selfdual Lagrangian on a reflexive Banach space $X$ such that  ${\rm Dom}_1( L)$ is closed and let $H_{L}$ be its Hamiltonian.  If  $\Lambda:D(\Lambda)\subset X\to X^{*}$  is a regular  map such that:
\begin{equation}
\label{one}
  {\rm Dom}_1( L)\subset D(\Lambda) \quad \hbox{\rm and $ \lim\limits_{\|x\|\to +\infty}H_{L}(0,-x)+\langle \Lambda x, x\rangle=+\infty,$}
 \end{equation}
 then the functional  $I(x)=L(x, \Lambda x) +\langle \Lambda x, x\rangle$ attains its minimum at $\bar x \in D(\Lambda)$ in such a way that:
  \begin{eqnarray}
 I( \bar x) &=&\inf_{x\in X} I(x)= 0\label{eqn:zero1}\\
  \hfill 0&\in&\Lambda \bar x + \partial L (\bar x).\label{eqn:zero2}
  \end{eqnarray}
\end{theorem}
We have denoted here the {\it effective domain} of $L$ by ${\rm Dom}(L)=\{(x,p)\in X\times X^*; L(x,p)<+\infty \}$, and by 
${\rm Dom}_1(L)$ its projection on $X$, that is  
 $
{\rm Dom}_1(L)=\{x\in X; L(x,p)<+\infty \hbox{ for some $p\in X^*$} \}.$
 
  The  Hamiltonian $H_L: X\times X \to \bar \R$ of $L$ is defined 
   by:
  \[
H_L(x,y)=\sup\{\langle y,p\rangle -L(x,p); p\in X^*\}, 
\]
which is the Legendre transform in the second variable. 

As shown in \cite{G3}, Theorem \ref{main.1}  applies readily to many nonlinear stationary equations giving variational proofs of existence of solutions. For example, one can obtain solutions of the incompressible stationary stationary Navier-Stokes equation on a smooth bounded domain $\Omega$ of $\R^{3}$
 \begin{equation}
\label{NSE1}
 \left\{ \begin{array}{lcl}
    \hfill
 (u\cdot \nabla)u +f &=&\nu \Delta u - \nabla  p \quad \, \, \hbox{\rm on $ \Omega$}\\
\hfill {\rm div} u&=&0 \quad \quad  \quad \quad   \quad   \hbox{\rm on  $\Omega$}\\
\hfill u &=&0 \quad \quad  \quad \quad   \quad \hbox{\rm on $\partial \Omega$}\\
\end{array}\right.
\end{equation}
where $\nu >0$ and $f\in L^{p}(\Omega;\R^{3})$, as follows. Letting 
\begin{equation}
\label{Phi}
\Phi (u)=\frac{\nu}{2} \int_{\Omega}\Sigma_{j,k=1}^{3}(\frac {\partial u_{j} } {\partial x_{k}})^{2}\, dx+\int_{\Omega}\Sigma_{j=1}^{3}f_{j}u_{j}
\end{equation}
be the convex continuous function on the space $X=\{u\in H^{1}_0(\Omega; {\bf R}^{3}); {\rm div} v=0\}$, and $\Phi^{*}$ be its Legendre transform  on $X^*$, 
Equation (\ref{NSE1}) can then be reformulated as
 \begin{equation}
\label{NSE2}
 \left\{ \begin{array}{lcl}
    \hfill
 \Lambda u &=& -\partial \Phi (u)=\nu \Delta u - f-\nabla  p \\
\hfill  u&\in& X,  \\
\end{array}\right. 
\end{equation}
where  $\Lambda: X \to X^{*}$ is the regular nonlinear operator defined as
\begin{equation}\label{Lambda}
\langle \Lambda u, v\rangle =\int_{\Omega}\Sigma_{j,k=1}^{3}u_{k}\frac {\partial u_{j} } {\partial x_{k}}v_{j}\, dx=\langle (u\cdot \nabla)u,v\rangle.
\end{equation}
Theorem \ref{main.1} then readily yields that  if $p>\frac{6}{5}$, then the infimum of the  functional 
 \begin{equation}
{I}(u)=\Phi (u)+\Phi^{*}(-(u\cdot \nabla)u ) 
\end{equation}
 on $X$ is equal to zero, and is attained at a solution of  (\ref{NSE1}).
 Theorem \ref{main.1} does not however cover the case of nonlinear evolutions such as the Navier-Stokes equations. This is because of the prohibitive coercivity condition (\ref{one}) that is not satisfied by the corresponding selfdual functional on the relevant path space. We shall therefore prove  a similar result under a more relaxed  coercivity condition that will allow us to prove a selfdual variational principle that is more appropriate to nonlinear evolution equations.

For that, we shall consider   an evolution triple $X \subset H \subset X^*$ where  $H$ is a Hilbert
space equiped with $\braket{}{}$ as scalar product, and where $X$ is a dense vector subspace of $H$, that is a reflexive Banach space once equipped with its own norm $\| \cdot \|$.  Let $[0,T]$ be a fixed real interval  and  consider for $p, q >1$, the  Banach space
$L^{p}_{X}$
as well as the  space ${\cal X}_{p,q}$ of all functions in $L^{p}_{X}$ such that $\dot{u} \in L^{q}_{X^*}$, equipped with the norm 
\[
     \|u\|_{{\cal X}_{p,q}} = \|u\|_{L^p_X} + \|\dot u\|_{L^q_{X^*}}.
     \]
      Let now $J$ be the duality map from $X$ to $X^*$,  i.e., for every $u \in X,$  $Ju$  is the element of the dual $X^*$ that is uniquely determined by the relation
\begin{eqnarray}
{\braket {Ju} u}=\|u\|_X^2 \text {  and  } \|Ju\|_{X^*}=\|u\|_X.
\end{eqnarray}
It is well-known that $J$  is one to one and onto $X^*$, while being monotone  and continuous from $X$ (with its strong topology) to $X^*$ equipped with its weak topology. We shall  need the following notion. 
 \begin{definition} \rm Let $L$ be  a time-dependent selfdual Lagrangian on $[0,T] \times X \times X^*$, and let $\Lambda:{\cal X}_{p,q}\to L^q_{X^{*}}$ be a given map. 
 Say that $L$ is {\it $\Lambda$-coercive} 
 if for any  sequence $\{x_n\}_{n=1}^{\infty} \subseteq {\cal X}_{p,q}$ we have 
 \[
\lim\limits_{\|x_n\|_{{\cal X}_{p,q}}\to +\infty}\int_0^T \big [ L(t, x_n(t),\dot x_n(t)+\Lambda x_n(t)+\frac{1}{n}\|x_n\|^{p-2}Jx_n(t))
+\langle x_n(t), \Lambda x_n(t) \rangle +\frac {1}{n} \|x_n(t)\|^p\big ] \,dt=+\infty.
\]
\end{definition}
Here is one useful corollary of the variational principle we establish in section 3 for nonlinear evolutions.

\begin{theorem}  \label{main.00} Let $X\subset H\subset X^*$  be an evolution triple where $X$ is a reflexive Banach space, and $H$ is a Hilbert space. For $p>1$ and $q=\frac{p}{p-1}$, assume that $\Lambda:{\cal X}_{p,q}\to  L^q_{X^{*}}$  is  a  regular map such that for some nondecreasing continuous real function $w$, and $0\leq k<1$, it satisfies
 \begin{equation} \label{Lambda-10}
\hbox{$\|\Lambda x\|_{L^q_{X^*}}\leq k\|\dot x\|_{L^q_{X^*}}+ w(\|x\|_{L^p_{X}})$ for every $x\in {\cal X}_{p,q}$,}
\end{equation}
and
\begin{eqnarray}\label{Lambda-20}
\hbox{$\big |\int_0^T \langle \Lambda  x(t),  x(t)\rangle \, dt  \big | \leq  w(\|x\|_{L^p_X})$ for every $x\in {\cal X}_{p,q}$}.
\end{eqnarray}
 Let $\ell$ be an anti-selfdual Lagrangian on $H \times H$ that is bounded below with  $0 \in {\rm Dom}(\ell)$, and  let $L$ be a time dependent  anti-selfdual Lagrangian on  $[0,T] \times X \times X^*$ that  is $\Lambda$-coercive and such  that for some $C>0$ and $r>1$, we have 
\begin{equation}\label{L.C}
\hbox{$\int_0^TL(t, u(t), 0) dt \leq C(1+\|u\|^r_{L^p_X})$ for every $u\in L^p_X$.}
\end{equation} 
 The following functional 
\begin{equation}
I(u)= \int_0^T \Big [  L  (t,  u(t),\dot { u}(t)+\Lambda  u(t)) +\langle \Lambda  u(t),  u(t) \rangle \Big ] \, dt +\ell (u(0)- u(T), \frac { u(T)+ u(0)}{2})
\end{equation}
then attains its minimum at $v \in {\cal X}_{p,q}$ in such a way that $I(v)=\inf_{u\in  {\cal X}_{p,q}}I(u)=0$ and
\begin{eqnarray}\label{principle3}
 \left \{ \begin{array}{lcl}
  \hfill -\Lambda v(t)- \dot {v}(t)& \in & \bar \partial L (t,v(t)), \\ \label{eqn:zero2}
  \hfill  -\frac{v(0)+v(T)}{2}
 & \in & \bar \partial \ell\big(v(0)-v(T)).
\end{array}\right.
  \end{eqnarray}
\end{theorem}
Now while the main Lagrangian $L$ is expected to be smooth and hence its subdifferential coincides with its gradient,  and the differential inclusion is often an equation, it is crucial that the boundary Lagrangian $\ell$ be allowed to be degenerate so as its subdifferential can cover the various boundary conditions discussed below.

  As a consequence of the above theorem, we provide a variational resolution to evolution equations involving  nonlinear operators such as the Navier-Stokes equation with various boundary conditions. Indeed, by considering  
  \begin{equation}
\label{T-NS.in}
 \left\{ \begin{array}{lcl}
    \hfill
 \frac {\partial u}{\partial t}+(u\cdot \nabla)u +f &=&\nu \Delta u - \nabla  p \quad \hbox{\rm on $ \Omega\subset \R^n$},\\
\hfill {\rm div} \, u&=&0 \quad \hbox{\rm on  $\Omega$},\\
\hfill u&=&0 \quad \hbox{\rm on $\partial \Omega$},\\
\end{array}\right.
\end{equation}
where   $f\in L^{2}_{X^*}([0,T])$,   $X=\{u\in H^{1}_0(\Omega; {\bf R}^{n}); {\rm div} v=0\}$, and $H=L^2(\Omega)$, 
  we can associate the nonlinear operator equation
 \begin{equation}
\label{T-NS-eq}
 \left\{ \begin{array}{lcl}
  \hfill \frac{\partial u}{\partial t}+ \Lambda u   &\in& -\partial \Phi (t,u) \\
 \hfill  \frac{ u(0)+ u(T)}{2} & \in& -\bar \partial \ell (u(0)-u(T)).
\end{array}\right.
\end{equation}
where $\ell$ is any anti-selfdual Lagrangian on $H \times H$, while  $\Phi$  and $\Lambda$ are defined in (\ref{Phi}) and (\ref{Lambda}) respectively.

Note that $\Lambda$ maps $X$ into its dual $X^*$ as long as the dimension $N\leq 4$. On the other hand, if we lift  $\Lambda$  to path space by defining $(\Lambda u)(t)=\Lambda (u (t))$, we have the following well-known results:
\begin{itemize}
\item If $N=2$, then $\Lambda$ is a regular operator from ${\cal X}_{2,2}[0,T]$ into  $L^2_{X^*}[0, T]$.
\item However, if $N=3$, we then have that $\Lambda$ is a regular operator from ${\cal X}_{2,2}[0,T]$ into  $L^{4/3}_{X^*}[0, T]$.
\end{itemize}
We therefore distinguish the two cases. 
\begin{corollary}\label{T-NS2} Assuming $N=2$, $f$ in $ L^{2}_{X^*}([0,T])$, and $\ell$ to be an anti-selfdual Lagrangian on $H \times H$ that is bounded from below, then the infimum of the functional
\[
I(u)=\int_0^T \big [\Phi( t, u(t))+\Phi^*(t, -\dot{u}(t)- (u\cdot \nabla)u(t)) \big ] \,dt + \ell (u(0)-u(T), \frac{u(0)+ u(T)}{2})
\]
on ${\cal X}_{2, 2}$ is zero and is attained at a  solution  $u$ of  (\ref{T-NS.in}) that satisfies the following time-boundary condition: 
\begin{equation} \label{br}
 -\frac{u(0)+u(T)}{2}
  \in  \bar \partial \ell\big(u(0)-u(T)).
  \end{equation}
   Moreover,  $u$ verifies the following ``energy identity":
\begin{equation}\label{T-NS2-in}
\hbox {$\|u(t)\|_{H}^2+2 \int_0^t \big [\Phi ( t, u(t))+\Phi^*(t, -\dot{u}(t)- (u\cdot \nabla)u(t)) \big ] \,dt = \|u(0)\|_{H}^2$ for every $t \in [0,T].$}
\end{equation}
In particular, with appropriate choices for the boundary Lagrangian $\ell$, the solution $u$ can be chosen to  verify either one of the following boundary conditions:
\begin{itemize}
\item an initial value problem:  $u(0) = u_0$ where $u_0$ is a given function in $X$. 
\item a periodic orbit  :  $u(0) = u(T)$,
\item an anti-periodic orbit  :  $u(0) =- u(T)$.
 \end{itemize}
\end{corollary}

However,  in the three dimensional case, we have to settle for the following result. 

\begin{corollary}\label{T-NS3} Assume $N=3$, $f$ in $ L^{2}_{X^*}([0, T])$, and consider  $\ell$ to be a selfdual Lagrangian on $H \times H$ that is now coercive in both variables.  Then,
there  exists $u \in {\cal X}_{2, {\frac {4}{3}}} $ such that
\[
I(u)=\int_0^T \big [\Phi( t, u(t))+\Phi^*(t, -\dot{u}(t)- (u\cdot \nabla)u(t)) \big ] \,dt + \ell (u(0)-u(T), \frac{u(0)+ u(T)}{2})\leq 0, 
\]
 and  $u$ is  a weak solution of  (\ref{T-NS.in}) that satisfies the time-boundary condition (\ref{br}).
  Moreover,  $u$ verifies the following ``energy inequality": 
 \begin{equation}\label{T-NS3-in}
\frac {\|u(T)\|_H^2}{2}+ \int_0^T \big [\Phi ( t, u(t))+\Phi^*(t, -\dot{u}(t)- (u\cdot \nabla)u(t)) \big ] \,dt  \leq \frac {\|u(0)\|_H^2}{2}.
\end{equation}
In particular, with appropriate choices for the boundary Lagrangian $\ell$, the solution $u$ will verify either one of the following boundary conditions:
\begin{itemize}
\item an initial value problem:  $u(0) = u_0$.
\item a periodicity condition of the form:  $u(0) = \alpha u(T)$,  for any given $\alpha$ with $-1< \alpha <1.$
 \end{itemize}
 
\end{corollary}
The above results are actually particular cases of a much more general nonlinear selfdual variational principles which applies to both the stationary and to the dynamic case. It will be stated and established in full generality in the next section.

\section{Basic properties of selfdual functionals} 

Consider the Hamiltonian $H=H_L$ associated to an ASD Lagrangian $L$ on $X\times X^*$. It is easy to check that   $H:X\times X\to \R\cup\{+\infty\}\cup\{-\infty\}$ then satisfies:
 \begin{itemize}
 \item  for each $y\in X$, the function $H_y: x\to -H(x,y)$ from $X$ to $\R\cup\{+\infty\}\cup\{-\infty\}$ is convex;
 \item  the function $x\to H(-y,-x)$ is the convex lower semi-continuous envelope of $H_y$.
 \end{itemize}
 
 It readily follows that for such a  Hamiltonian, the function  $y\to H(x,y)$ is convex and lower semi-continuous for each $x\in X$, and that  the following inequality holds: 
 \begin{equation}
\label{almost.odd}
\hbox{$H(-y,-x)\leq -H(x,y)$ for every $(x,y)\in X\times X$.}
\end{equation}
In particular, we have
 \begin{equation}
\label{negative.diagonal}
\hbox{$H(x,-x)\leq 0$  for every $x\in X$.}
\end{equation}
Note that $H_L$ is always concave in the first variable, however, it is not necessarily upper semi-continuous in the first variable. 

Another property of ASD Lagrangians that will be used in the sequel is the fact that 
\begin{equation}\label{duality}
L(x,p)=(H_L)_2^*(x,p)=(-H_L)^*_1(-p,-x)
\end{equation}
where $(f)_1^*$ (resp., $(f)_2^*$) denotes the Legendre transform of  a function $f$ on $X\times X$, with respect to the first (resp., second) variable. It then follows that if we define the following operation on two ASD Lagrangians $L$ and $M$ on $X\times X^*$, 
\begin{equation}
L\oplus M (x,p)=\inf \{L(x, r)+M(x, p-r); r\in X^*\},
\end{equation}
then we have for any $(x,p)\in X\times X^*$, 
\begin{equation}\label{duality.bis}
L\oplus M (x,p)=\sup\{ \langle y, -p\rangle +H_L(y, -x)+H_M(y,-x); y\in X\}.
\end{equation} 
   As   in \cite{G3}, we consider the following notion which extends considerably the class of  Hamiltonians associated to selfdual Lagrangians.
\begin{definition} \rm Let $E$ be a convex  subset of a reflexive Banach space $X$. 
\begin{enumerate}
\item A functional $M: E\times E \to \R$ is said to be {\it an  anti-symmetric  Hamiltonian} on $E\times E$ if it satisfies the following conditions:
\begin{equation}
\hbox{For every $x\in E$, the function $y\to M(x, y)$ is concave on $E$.}
\end{equation}
\begin{equation}
\hbox{$M(x,x) \leq 0$ for every $x\in E$.}
\end{equation}
\item It is said to be a {\it regular  anti-symmetric Hamiltonian} if in addition it satisfies:
\begin{equation}
\hbox{For every $y\in E$, the function $x\to M(x,y)$ is weakly lower semi-continuous on $E$.}
   \end{equation}
 \end{enumerate}
  \end{definition}
   The class of {\it regular anti-symmetric Hamiltonians} on a given convex set $E$ --denoted  ${\cal H}^{asym}(E)$-- is an interesting class of its own. It  contains the ``Maxwellian" Hamiltonians  $H(x,y)=\phi (y) -\phi (-x)+\langle Ay, x\rangle$, where $\phi$ is convex and $A$ is skew-adjoint. 
More generally, 
\begin{enumerate}
\item  If $L$ is an anti-selfdual Lagrangian on a Banach space $X$, then the Hamiltonian $M(x,y)=H_L(y,-x)$ is in ${\cal H}^{asym}(X)$. 

\item If  $\Lambda: D(\Lambda)\subset X\to X^*$ is a --non necessarily linear-- regular, 
 then the Hamiltonian $H(x,y)=\langle x-y, \Lambda x\rangle $ is in ${\cal H}^{asym}(D(\Lambda))$. 
\end{enumerate}

Since ${\cal H}^{asym}(X)$ is obviously a convex cone, we can therefore superpose certain non-linear operators with anti-selfdual Lagrangians, via their corresponding anti-symmetric Hamiltonians, to obtain a remarkably rich family that generates non-convex selfdual functionals as follows. 

\begin{definition} A functional  $I: X\to \R \cup \{+\infty\}$  is said to be {\it selfdual on a convex set $E\subset X$} if it is non-negative and if there exists a regular anti-symmetric  Hamiltonian $M: E\times E \to \R$ such that for every $x\in E$, 
\begin{equation}
\hbox{$I(x)=\sup\limits_{y\in E}M(x,y)$.}
\end{equation}
\end{definition}
A key aspect of our variational approach is that solutions of many nonlinear PDEs can be obtained by minimizing  properly chosen selfdual functionals in such a way that the infimum is actually zero.
This is indeed the case in view of the following immediate application of a fundamental min-max theorem of Ky-Fan (see \cite{G3}).  

\begin{proposition} \label{SD1} Let  $I: E\to \R \cup \{+\infty\}$ be  a selfdual functional  on a closed convex  subset $E$ of a reflexive Banach space $X$, with $M$ being its corresponding anti-symmetric Hamiltonian on $E\times E$. If $M$ is coercive in the following sense
\begin{equation}
\lim_{\|x\|\to +\infty}M(x, 0)=+\infty, 
\end{equation}
then there exists $\bar x\in E$ such that 
$
I(\bar x)=\sup\limits_{y\in E}M(\bar x,y)= 0.
$
\end{proposition}

The following was also proved in \cite{G3}.
  \begin{proposition} \label{lift} Let $X\subset H \subset X^*$ be an evolution pair and consider 
a time-dependent  anti-selfdual Lagrangian $L$ on $[0,T]\times X\times X^*$ such that  
\begin{equation}
\label{L.A}
  \hbox{   For each $p\in L^{q}_{X^*}$, the map $u\to \int_0^T L(t, u(t), p (t)) dt$ is continuous on $L^{p}_{X}$ }
 \end{equation}
 \begin{equation}
 \label{L.B}
  \hbox{   The map $u\to \int_0^T L(t, u(t), 0) dt$ is bounded on the unit ball of $L^{p}_{X}$. }
 \end{equation}
Let $\ell$ be an anti-selfdual  Lagrangian on $H\times H$ such that:
  \begin{equation} \label{ell.A}
  \hbox{ $ -C\leq  \ell (a,b) \leq C(1+\|a\|_H^2+\|b\|^2_H)$ for all $(a,b)\in H\times H$.}
   \end{equation}
 Then the Lagrangian
\[
{\cal L}(u,p)= \left\{ \begin{array}{lcl}
  \hfill  \int_0^T L(t, u(t), p (t)+\dot{u}(t)) dt + \ell (u(0)-u(T), \frac{u(T)+u(0)}{2})\, &{\rm if}& u\in {\cal X}_{p,q}\\
  +\infty &\,& {\rm otherwise} \\
\end{array}\right.
\]
is anti-selfdual on $L^{p}_{X}\times L^q_{X^*}$.
\end{proposition}

Consider now the following convex lower semi-continuous function on $L^p_X$:
 \begin{equation}
\label{dot}
 \psi(u)=\left\{ \begin{array}{lcl}
 \frac {1}{q}\int_0^T\|\dot u (t)\|^q_{X^*}\, dt
  &{\rm if }& u\in  {\cal X}_{p,q} \\
 +\infty & {\rm if } & u\in L^p_X\setminus  {\cal X}_{p,q}, 
\end{array}\right.
\end{equation}
and for any $\mu>0$, we let $\Psi_\mu$ be the anti-selfdual Lagrangian on $L^{p}_{X}\times L^{q}_{X^*}$ defined by 
 \begin{equation}
\Psi_\mu(u, r)= \mu \psi (u)+ {\mu} \psi^* (-\frac{r}{\mu}). 
\end{equation}

Now for each  $(u,r) \in  L^p_X  \times L^q_{X^*},$ define
\begin{eqnarray}\label{Phi-2}
{\cal L} \oplus  \Psi_{\mu} (u,r):=\inf_{s \in L^q_{X^*} } \big \{{\cal L} (u,s )+  \Psi_{\mu} (u,r-s) \big \}
\end{eqnarray}
\begin{lemma} \label{double.dot} Let $L$ and $\ell$ be two anti-selfdual Lagrangians verifying the hypothesis of Proposition \ref{lift}, and let ${\cal L}$ be the corresponding anti-selfdual Lagrangian on path space  $L^p_X\times L^p_X$.  Suppose $\Gamma$ is a regular operator from ${\cal X}_{p,q}$ into $L^q_{X^*}$ then,  
\begin{enumerate}
\item The functional 
\[
I_\mu (u)=  {\cal L} \oplus  \Psi_{\mu}
(u,\Gamma u)+\int_0^T \langle \Gamma u (t), u  (t)\rangle  \, dt 
\]
is selfdual on  ${\cal X}_{p,q} \times {\cal X}_{p,q}$, and its corresponding anti-symmetric Hamiltonian on ${\cal X}_{p,q} \times {\cal X}_{p,q}$ is  
\begin{eqnarray*}
M_\mu(u,v):=&\int_0^T \langle \Gamma u (t), u(t)-v (t) \rangle dt+H_{\cal L}(v, -u)
+\mu \psi (u)-\mu \psi (v), 
\end{eqnarray*}
where $H_{\cal L}(v, u)=\sup_{p \in L^q_{X^*}} \{ \int_0^T \langle p,u \rangle \, dt -{\cal L} (v,p)$ is the Hamiltonian of ${\cal L}$ on $L^p_X\times L^p_X$. \\

\item  If in addition 
$
\lim\limits_{\|u\|_{{\cal X}_{p,q}}\to +\infty}\int_0^T \langle \Gamma u (t), u(t)\rangle dt+H_{\cal L}(0, -u)
+\mu \psi (u)=+\infty, 
$
then there exists $u \in  {\cal X}_{p,q}$ with  $\partial \psi (u) \in L^q_{X^*}$  such that
\begin{eqnarray} \label{app-1}
\dot u (t)+\Gamma u (t)+ \mu  \partial \psi (u (t))&\in & -\bar \partial   L(t,u (t))\\ 
\frac {u(T)+u (0)}{2}&\in &-\bar \partial \ell (u (0)-u (T))\\  
 \dot u (T)&=& \dot u (0)=0.
\end{eqnarray}
\end{enumerate}
\end{lemma} 
\noindent {\bf Proof:}  First note that since ${\cal L}$ and $\Psi_\mu$ are anti-selfdual, it is easy to see that 
${\cal L} \oplus \Psi_{\mu} (u, p) +\langle u, p\rangle\geq 0$ for all $(u, p) \in  L^p_X\times L^q_{X^*}$, and therefore $I (u) \geq 0$ on ${\cal X}_{p,q}$.  

Now by   (\ref{duality.bis}), we have for any  $(u, p) \in  L^p_X\times L^q_{X^*}$, 
\begin{eqnarray*}
{\cal L} \oplus \Psi_{\mu} (u, p)&=&\sup_{v \in L^p_{X}} \{ \int_0^T \langle -p,v \rangle \, dt +H_{\cal L} (v,-u)+\mu \psi (-u)-\mu \psi (v) \}.
\end{eqnarray*}
But for $u\in {\cal X}_{p,q}$ and $v \in L^p_X \setminus {\cal X}_{p,q}$, we have 
$H_{\cal L}(v, -u)=\sup_{p \in L^q_{X^*}} \{ \int_0^T -\langle p,u \rangle \, dt - {\cal L} (v,p) \}=-\infty,$
and therefore for any $u\in {\cal X}_{p,q}$, we have
\begin{eqnarray*}
\sup_{v \in {\cal X}_{p,q}}M_\mu(u,v)&=&\sup_{v \in L^p_X}M_\mu(u,v)\\
&=&\int_0^T \langle \Gamma u (t), u(t) \rangle dt+\sup_{v \in L^p_X}\int_0^T \langle \Gamma u (t), -v (t) \rangle dt+H_{\cal L}(v, -u)
+\mu \psi (u)-\mu \psi (v)\\
&=&
\int_0^T \langle \Gamma u (t), u(t) \rangle dt+{\cal L} \oplus \Psi_{\mu} (u, \Gamma u)\\
&=&
I(u).
\end{eqnarray*}
It follows from Proposition \ref{SD1} that there exists $u_\mu \in {\cal X}_{p,q}$ such that  
\begin{eqnarray}
I_\mu (u_\mu )=  {\cal L} \oplus  \Psi_{\mu}
(u_\mu ,\Gamma u_\mu )+\int_0^T \langle\Gamma  u_\mu  (t), u_\mu   (t)\rangle  \, dt=0. 
\end{eqnarray}
Since 
 ${\cal L} \oplus  \Psi_{\mu} (u,p)$
 is convex and coercive in the second variable, there exists $p \in  L^q_{X^*}$ such that 
\begin{eqnarray}
{\cal L} \oplus  \Psi_{\mu}
(u_\mu ,\Gamma u_\mu )={\cal M}_{\cal L} (u,p)+\Psi_{\mu}(u, \Gamma u_\mu-p).
\end{eqnarray}
It follows that 
\begin{eqnarray*}
0 & =&  
{\cal L} (u_\mu,p)+\Psi_{\mu}(u_\mu, \Gamma u_\mu-p)+\int_0^T \langle \Gamma u_\mu(t), u_\mu(t) \rangle  \, dt \\ &= & \int_0^T \big [ L (t, u_\mu(t),\dot u_\mu (t)+p(t))+ \langle u_\mu(t),p(t) \rangle  \big ]\,dt+ \ell(u_\mu(T)-u_\mu(0), \frac{u_\mu(T)+u_\mu(0)}{2})\\ && +\Psi_{\mu}(u_\mu, \Gamma u_\mu-p)+\int_0^T \langle \Gamma u_\mu (t)-p(t), u_\mu   (t)\rangle  \, dt\\  & = & \int_0^T \big [ L (t, u_\mu(t),\dot u_\mu (t)+p(t))+ \langle u_\mu(t), \dot u_\mu(t)+ p(t) \rangle  \big ]\,dt\\
&&-\frac{1}{2}\|u_\mu(T)\|^2+\frac{1}{2}\|u_\mu(0)\|^2+\ell(u_\mu(0)-u_\mu(T), \frac{u_\mu(T)+u_\mu(0)}{2})\\ &&  +\Psi_{\mu}(u_\mu, \Gamma u_\mu-p)+\int_0^T \langle \Gamma u_\mu-p, u_\mu   (t)\rangle  \, dt.
\end{eqnarray*}
Since this is the sum of three non-negative terms, we get the following three identities, 
\begin{equation}
\hbox{$\int_0^T \big [ L (t, u_\mu(t),\dot u_\mu (t)+p(t))+ \langle u_\mu, \dot u_\mu+ p \rangle  \big ]\,dt=0$, }
\end{equation}
\begin{equation}
\hbox{$\Psi_{\mu}(u_\mu, \Gamma u_\mu-p)+\int_0^T \langle \Gamma u_\mu-p, u_\mu   (t)\rangle  \, dt=0$, }
\end{equation}
\begin{equation}
 \ell(u_\mu(0)-u_\mu(T), \frac{u_\mu(T)+u_\mu(0)}{2})-\frac{1}{2}\|u_\mu(T)\|^2+\frac{1}{2}\|u_\mu(0)\|^2=0.
\end{equation}
It follows from the limiting case of Fenchel duality that
\begin{eqnarray*} \label{app-1}
 \dot u_\mu (t)+ \Gamma u_\mu (t)+ \mu  \partial \psi (u_\mu (t)))&\in & -\bar \partial   L(t,u_\mu (t)) \hbox{for almost all $t\in [0, T]$}\\
 \frac{u_\mu(T)+u_\mu(0)}{2}&\in &-\bar \partial \ell (u_\mu(0)-u_\mu(T)).  
\end{eqnarray*}
Since  $u:=u_\mu \in  {\cal X}_{p,q}$, we have  that $-\mu  \partial \psi (u (t)))=\dot u (t)+\Gamma u (t)+ \bar \partial   L(t,u (t)) \in L^q_{X^*}$.

It follows that  $\partial \psi (u (t)))=- \frac{d}{dt}(\|\dot u\|_*^{q-2}J_T^{-1}\dot u)$, where $J_T$ is the duality map between $L^p_X$ and $L^q_{X^*}$. Hence, for each $v \in {\cal X}_{p,q} $ we have 
\begin{eqnarray*} \label{app-1}
0&=& \int_0^T\Big [ \langle \dot u (t)+\Gamma  u (t)+\bar \partial   L(t,u (t)), v \rangle +\mu  \langle \|\dot u\|_*^{q-2}J_T^{-1}\dot u, \dot v \rangle  \Big ]\, dt \\ &=& \int_0^T \langle \dot u (t)+\Gamma u (t)- \mu \frac{d}{dt}(\|\dot u\|_*^{q-2}J_T^{-1}\dot u)+\bar \partial   L(t,u (t)), v \rangle \, dt \\ &&+\mu \langle \|\dot u(T)\|_*^{q-2}J_T^{-1}\dot u(T), v(T) \rangle -\mu \langle \|\dot u(0)\|_*^{q-2}J_T^{-1}\dot u(0), v(0) \rangle
\end{eqnarray*}
from which we deduce that 
\begin{eqnarray*}
 \dot u (t)+\Gamma u (t)- \frac{d}{dt}(\|\dot u\|^{q-2}J_T^{-1}\dot u(t))&\in &-\bar \partial   L(t,u (t))\\
 \dot u(T)&=&\dot u(0)=0.
 \end{eqnarray*} 
\hfill $\square$

We shall  make repeated use of the following  lemma which describes three ways of regularizing an anti-selfdual Lagrangian by way of $\lambda$-convolution. It is an immediate consequence of the calculus of anti-selfdual Lagrangians developed in \cite{G2} to which we refer the reader. 

 \begin{lemma} \label{2-reg} For a  Lagrangian  $L: X \times X^*\rightarrow \mathbb{R}\cup \{+\infty\}$, define for every $(x, r)\in X\times X^*$
\begin{eqnarray*}
L^1_{\lambda}(x,r) =\inf \{ L(y,r)+ \frac {\|x-y\|^p}{\lambda p}+ \frac {\lambda ^{q-1}\|r\|^q}{ q}; y \in X \}
\end{eqnarray*}
and
\begin{eqnarray*}
L^2_{\lambda}(x,r) =\inf \{ L(x,s)+ \frac {\|r-s\|^q}{\lambda q}+ \frac {\lambda ^{p-1}\|x\|^p}{p}; s \in X^* \}
\end{eqnarray*}
and 
\begin{eqnarray*}
L^{1,2}_{\lambda}(x,r)=\inf \big \{   L(y,s)+ \frac{1}{2 \la}\|x-y\|^2+  \frac{\la}{2 }\|r\|_*^2+ \frac{1}{2 \la}\|s-r\|_*^2+  \frac{\la}{2 }\|y\|^2; \, y \in X, s \in X^*\big \}
\end{eqnarray*}

If $L$ is anti-selfdual then the following hold:
\begin{enumerate}
\item   $L^1_{\lambda}$, $ L^2_{\lambda}$ and $L^{1,2}_{\lambda}$ are  also anti-selfdual Lagrangians on $X \times X^*$. 
\item  $L^1_{\lambda}$ (resp., $L^2_{\lambda}$) (resp., $L^{1,2}_{\lambda}$) is continuous in the first variable (resp., in the second variable) (resp., in both variables). Moreover,  $\|\bar \partial L^1_{\lambda}( x) \| \leq \frac {\|x\|}{\lambda}$ for every $x \in X $.
\item  $\bar  \partial L^2_{\lambda} (x)=\bar  \partial L (x)+\lambda^{p-1} \|x\|^{p-2}Jx$ for every $x \in X $.
\item   $\bar  \partial L^1_{\lambda} (x)= \bar  \partial L (x+\lambda^{q-1} \|r\|^{q-2}J^{-1}r )$ for every $x \in X $ where $r= \bar \partial L(x)$.
\item Suppose $L$ is bounded from below.  If $x_{\la}\rightharpoonup x$ and $p_{\la}\rightharpoonup p$ weakly in $X$ and $X^*$ respectively as $\lambda \to 0$,  and if  $L^{1,2}_{\lambda}(x_{\la},p_{\la})$ is bounded from above,  then
\begin{eqnarray*}
L(x,p)\leq \liminf_{\la\to 0}L^{1,2}_{\lambda}(x_{\la},p_{\la}).
\end{eqnarray*}

\end{enumerate}
\end{lemma}
{\bf Proof:} It suffices to notice that  $L_\lambda^1=L\star M_\lambda$ and $L_\lambda^2=L\oplus M_\lambda$  where $M_\lambda(x,r)=\psi_\lambda (x)+\psi_\lambda^*(r)$ with $\psi_\lambda(x)=\frac{1}{\lambda p}\|x\|^p$. Note that 
$L^{1,2}_\lambda= (L\oplus M_\lambda)\star M_\lambda$ with $M_\lambda (x,r)= \frac{1}{2\lambda}\|x\|^2 +\frac{\lambda}{2}\|r\|^2$. The rest follows from the calculus of selfdual Lagrangians developed in \cite{G2}.  \hfill $\square$\\

\section{A selfdual variational principle for nonlinear evolutions}
 
This section is dedicated to the proof of the following  general variational principle for nonlinear evolutions.  

\begin{theorem}  \label{main.20} Let $X\subset H\subset X^*$  be an evolution triple where $X$ is a reflexive  Banach space, and $H$ is a Hilbert space. Let $L$ be a time dependent  anti-selfdual Lagrangian on  $[0,T] \times X \times X^*$ such that for some $C>0$ and $r>0$, we have 
\begin{equation}\label{L.C}
\hbox{$\int_0^TL(t, u(t), 0) dt \leq C(1+\|u\|^r_{L^p_X})$ for every $u\in L^p_X$.}
\end{equation} 
 Let $\ell$ be an anti-selfdual Lagrangian on $H \times H$ that is bounded below with  $0 \in {\rm Dom}(\ell)$, and consider 
 $\Lambda:{\cal X}_{p,q}\to  L^q_{X^{*}}$  to be  a  regular map such that :
\begin{equation} \label{Lambda-1}
\hbox{$\|\Lambda u\|_{L^q_{X^*}}\leq k\|\dot u\|_{L^q_{X^*}}+ w(\|u\|_{L^p_{X}})$ for every $u\in {\cal X}_{p,q}$,}
\end{equation}
where  $w$ is a nondecreasing continuous real function and $0<k<1$. Assume that $L$ is $\Lambda$-coercive and 
that one of the following two conditions hold:
\begin{itemize}
\item[(A)] 
$\big |\int_0^T \langle \Lambda  u(t),  u(t)\rangle \, dt  \big | \leq  w(\|u\|_{L^p_X})$ for every $u\in {\cal X}_{p,q}$.
\item[(B)] For each $p \in L^q_{X^*}$, the functional $u\rightarrow \int_0^T L(t, u(t), p(t)) \, dt$ is continuous on $L^p_X$, and there exists  $C>0$   such that for every $u\in L^p_X$ we have:
\begin{eqnarray}
\|\bar \partial L(t,u)\|_{L^q_{X^*}} &\leq&  w(\|u\|_{L^p_X}), \label{partial L}\\
\int_0^T\langle  \bar \partial L(t,u(t))+\Lambda  u(t),  u(t)\rangle  \, dt  &\geq& -C(\|u\|_{L^p_X}+1).\label{coercive-2}
\end{eqnarray}
\end{itemize}
Then  the functional 
\begin{equation}
I(u)= \int_0^T \Big [  L  (t,  u(t),\dot { u}(t)+\Lambda  u(t)) +\langle \Lambda  u(t),  u(t) \rangle \Big ] \, dt +\ell (u(0)- u(T), \frac { u(T)+ u(0)}{2})
\end{equation}
attains its minimum at $v \in {\cal X}_{p,q}$ in such a way that $I(v)=\inf_{u\in  {\cal X}_{p,q}}I(u)=0$ and
\begin{eqnarray}\label{principle3}
 \left \{ \begin{array}{lcl}
  \hfill -\Lambda v(t)- \dot {v}(t)& = & \bar \partial L (t,v(t)), \\ \label{eqn:zero2}
  \hfill  -\frac{v(0)+v(T)}{2}
 & \in & \bar \partial \ell\big(v(0)-v(T)).
\end{array}\right.
  \end{eqnarray}
\end{theorem}
For the proof, we start with the following proposition in  which we consider a regularization (coercivization) of the anti-selfdual Lagrangian ${\cal L}$ by the ASD Lagrangian $\Psi_\mu$, and also a  perturbation of the operator $\Lambda$ by  operator  
\begin{equation}
 K u=w(\|u\|_{L^p_X})J_Tu+\|u\|_{L^p_X}^{p-1} J_Tu
 \end{equation}
  which is regular from ${\cal X}_{p,q}$ into $L^q_{X^*}$. 

\begin{lemma} \label{app-lem}  Let $\Lambda$ be a regular  map from ${\cal X}_{p,q}$ into $L^q_{X^*}$ satisfing (\ref{Lambda-1}). Let $L$ to be a time-dependent anti-selfdual Lagrangian on $[0,T]\times X\times X^*$, satisfying conditions (\ref{L.A}) and  (\ref{L.B}) and  let $\ell$ be  an  anti-selfdual Lagrangian on $H \times H$ satisfying condition  (\ref{ell.A}).  
Then 
for any $\mu >0$, the functional 
\[
I_\mu (u)=  {\cal L} \oplus  \Psi_{\mu}
(u,\Lambda u+K u)+\int_0^T \langle \Lambda u (t)+ K u (t), u  (t)\rangle  \, dt 
\]
is selfdual on  ${\cal X}_{p,q} \times {\cal X}_{p,q}$. 
Moreover,  there exists $u_{\mu} \in \{u \in {\cal X}_{p,q};  \partial \psi (u) \in L^q_{X^*}, \dot u(T)= \dot u(0)=0\}$ such that
\begin{eqnarray} \label{app-1}
\dot u_{\mu} (t)+\Lambda u_{\mu} (t)+ K u_{\mu} (t)+ \mu  \partial \psi (u_{\mu} (t)))&\in & -\bar \partial   L(t,u_{\mu} (t))\\
 \ell (u_{\mu}(0)-u_{\mu}(T),  \frac {u_{\mu}(T)+u_{\mu}(0)}{2})&=&\int_0^T \langle \dot u_{\mu} (t) ,u_{\mu}  (t)\rangle \, dt
\end{eqnarray}
\end{lemma}
\noindent{\bf Proof:} It suffices to apply Lemma \ref{double.dot} to the regular operator $\Gamma =\Lambda +K$, provided we show the required coercivity condition $\lim\limits_{\|u\|_{{\cal X}_{p,q}}\to +\infty}\ M(u,0)=+\infty$ where
\[
M(u,0)=\int_0^T \langle \Lambda u (t) +K u(t), u(t)\rangle dt+H_{\cal L}(0, -u)
+\mu \psi (u).
\]
Note first that it follows from (\ref{Lambda-1})  that for $\epsilon < \frac{\mu}{q}$, there exists $C(\epsilon)$ such that 
\begin{eqnarray*}
\int_0^T \langle \Lambda u (t), u (t)\rangle \,dt &\leq  &k \|u\|_{L^p_X}\|\dot u\|_{L^q_{X^*}}+w(\|u\|_{L^p_X}) \|u\|_{L^p_X}\\
&\leq & \epsilon \|\dot u\|^q_{L^q_{X^*}}+C(\epsilon)\|u\|^p_{L^p_X}+w(\|u\|_{L^p_X}) \|u\|_{L^p_X}.
\end{eqnarray*}
On the other hand, by the definition of $K$, we have
\begin{eqnarray*}
\int_0^T \langle K u (t), u (t)\rangle \,dt =w(\|u\|_{L^p_X}) \|u\|_{L^p_X}^2+ \|u\|_{L^p_X}^{p+1}.
\end{eqnarray*}
Therefore the coercivity follows from the following estimate:
\begin{eqnarray*} \label{app-4}
M(u,0)&=&\int_0^T \Big [\langle \Lambda u (t)+ K u (t) , u (t) \rangle \Big] \, dt+ H_{\cal L}(u ,0)+\mu \frac {1}{q}\|\dot u\|_{L^q_{X^*}}^ q\\ 
&\geq &-  \epsilon \|\dot u\|^q_{L^q_{X^*}}-C(\epsilon)\|u\|^p_{L^p_X}-w(\|u\|_{L^p_X}) \|u\|_{L^p_X}+w(\|u\|_{L^p_X}) \|u\|_{L^p_X}^2+ \|u\|_{L^p_X}^{p+1}\\ &&- H_{\cal L}(0 ,0)+\mu \frac {1}{q}\|\dot u\|_{L^q_{X^*}}^ q\\
&\geq &(\frac{\mu}{q}-  \epsilon) \|\dot u\|^q_{L^q_{X^*}}+ \|u\|^{p+1}_{L^p_X}\big(1+ o(\|u\|_{L^p_X})\big).
\end{eqnarray*}

In the following lemma, we get rid of the regularizing diffusive term $\mu \psi (u)$ and prove the theorem with $\Lambda$ replaced by the operator $\Lambda +K$, and  under the additional assumption that $\ell$ satisfies the boundedness condition  (\ref{ell.A}).  

\begin{lemma}\label{pert-prop}  Let $L$ be  a time dependent anti-selfdual Lagrangian as  in Theorem  \ref{main.20} satisfying  either one of  conditions $(A)$  or $(B)$, and assume  that  $\ell $ is  an anti-selfdual Lagrangian on $H \times H$ that satisfies condition  (\ref{ell.A}).
Then there exists $u \in {\cal X}_{p,q} $ such that
\begin{eqnarray*}
\int_0^T\big [ L(t,u(t), \dot u(t)+ \Lambda u(t)+ K u(t)) \, dt+\langle \Lambda u(t)+Ku(t),
u(t)\rangle \big ] \, dt+\ell (u(0)- u(T), \frac { u(T)+ u(0)}{2})=0.
\end{eqnarray*}
\end{lemma}

\noindent{\bf Proof  under condition (B):} Note first that in this case $L$ satisfies both conditions (\ref{L.A}) and (\ref{L.B})  of Lemma  \ref{double.dot}, which then yields  for every $\mu>0$ an element  $u_{\mu} \in {\cal X}_{p,q}$  satisfying
\begin{eqnarray} \label{cdn-B-1}
\dot u_{\mu} (t)+\Lambda u_{\mu} (t)+ K u_{\mu} (t)+ \mu  \partial \psi (u_{\mu} (t)))\in  -\bar \partial   L(t,u_{\mu} (t)), 
\end{eqnarray}
and
\begin{eqnarray}\label{cdn-B-2}
\ell ( u_\mu(0)-u_\mu(T), \frac {u_\mu(T)+u_\mu(0)}{2})&=&\int_0^T \langle \dot u_\mu (t) ,u_\mu  (t)\rangle \, dt.
\end{eqnarray}
We now establish upper bounds on the norm of $u_{\mu}$ in ${\cal X}_{p,q}$.  Multiplying (\ref{cdn-B-1}) by $u_{\mu}$ and
 integrating over $[0,T]$  we obtain
\begin{eqnarray}\label{cdn-B-3}
\int _0^T \langle \dot u_{\mu} (t)+\Lambda u_{\mu} (t)+ K u_{\mu} (t)+ \mu  \partial \psi (u_{\mu} (t)),u_{\mu} (t) \rangle \, dt  =- \int _0^T \langle \bar \partial   L(t,u_{\mu} (t)), u_{\mu} (t) \rangle \, dt.
\end{eqnarray}
It follows from (\ref{coercive-2}) and the above equality that
\begin{eqnarray}\label{cdn-B-30}
 \int_0^T
\langle \dot u_{\mu} (t)+ K u_{\mu}  (t)+\mu\partial \psi (u_{\mu} (t)),u_{\mu} (t) \rangle  \leq C(1+\|u_{\mu}\|_{L^p_X}).
\end{eqnarray}
Taking into account (\ref{cdn-B-2}) and the fact that $\int_0^T\partial \psi (u_{\mu} (t)),u_{\mu} (t) \rangle \geq 0$, it follows from (\ref{cdn-B-30}) that
\begin{eqnarray*}
\ell ( u_{\mu}(0)-u_{\mu}(T), \frac {u_{\mu}(T)+u_{\mu}(0)}{2})+ \int_0^T
\langle  K u_{\mu}  (t),u_{\mu} (t) \rangle  \leq C(1+\|u_{\mu}\|_{L^p_X}).
\end{eqnarray*}
Since $\ell$ is bounded from below (say by $C_1$), the above inequality implies that  $\|u_{\mu}\|_{L^p_X}$
is bounded, since then we have 
\begin{eqnarray*}
C_1+ w(\|u_\mu\|_{L^p_X}) \|u_\mu\|_{L^p_X}^2+ \|u_\mu\|_{L^p_X}^{p+1}  \leq C\|u_{\mu}\|_{L^p_X}.
\end{eqnarray*}

 Now we show that $\|\dot u_{\mu}\|_{L^q_{X^*}}$ is also bounded.
For that, we  multiply (\ref{cdn-B-1}) by $J_{T}^{-1} \dot u_{\mu}$ to get that
\begin{eqnarray}\label{cdn-B-4}
\|\dot u_{\mu}\|^2_{L^q_{X^*}}+\int_0^T \big [ \langle \Lambda u_{\mu}(t)+
K u_{\mu}(t)+ \mu  \partial \psi (u_{\mu}(t))+\bar \partial   L(t,u_{\mu}(t)), J^{-1}  \dot u_{\mu}(t) \rangle \big ]\,dt   =0.
\end{eqnarray}
The last identity and the fact that $\int_0^T\langle \partial \psi (u_\mu (t)), J^{-1}  \dot u_{\mu}(t) \rangle \,dt =0$ imply 
that
\begin{eqnarray*}
\|\dot u_{\mu}\|^2_{L^q_{X^*}}\leq  \|\Lambda u_{\mu}\|_{L^q_{X^*}}
\| \dot u_{\mu}\|_{L^q_{X^*}}+ \|K u_{\mu}\|_{L^q_{X^*}} \| \dot u_{\mu}\|_{L^q_{X^*}}+C\| \dot u_{\mu}\|_{L^q_{X^*}}.
\end{eqnarray*}
It follows from the above inequality and (\ref{Lambda-1}) that
\begin{eqnarray*}
\|\dot u_{\mu}\|_{L^q_{X^*}}\leq  \|\Lambda u_{\mu}\|_{L^q_{X^*}}+
\|K u_{\mu}\|_{L^q_{X^*}} +C \leq k\|\dot u_{\mu}\|_{L^q_{X^*}}+ w(\|u\|_{L^p_X})+\|K u_{\mu}\|_{L^q_{X^*}}
\end{eqnarray*}
from which we obtain that
$(1-k)\|\dot u_{\mu}\|_{L^q_{X^*}}\leq  w(\|u_{\mu}\|_{L^p_X})+\|K u_{\mu}\|_{L^q_{X^*}}$, 
which means that $\|\dot u_{\mu}\|_{L^q_{X^*}}$ is bounded. 

Consider now $u\in {\cal X}_{p,q}$ such that $u_{\mu} \rightharpoonup u$
weakly in ${L^p_{X}}$ and  $\dot {u}_{\mu} \rightharpoonup \dot u$ in ${L^q_{X^*}}$. From (\ref{cdn-B-1}) and (\ref{cdn-B-2}) we have
\begin{eqnarray*}
J_\mu (u_\mu):&=&\int_0^T \Big [\langle \Lambda u_{\mu}(t)+ K u_{\mu}(t), u_{\mu}(t) \rangle + L (t,u_{\mu}(t),\dot u_{\mu}(t)+
 \Lambda u_{\mu}(t)+ K u_{\mu}(t)+ \mu  \partial \psi (u_{\mu}(t)) )\Big] \, dt\\
 &&+\ell (u_\mu(0)- u_\mu(T), \frac { u_\mu(T)+ u_\mu(0)}{2})\\
&\leq &\int_0^T \Big [\langle \Lambda u_{\mu}(t)+ K u_{\mu}(t)+\mu  \partial \psi (u_{\mu}(t)), u_{\mu}(t) \rangle + L (t,u_{\mu}(t),\dot u_{\mu}(t)+
 \Lambda u_{\mu}(t)+ K u_{\mu}(t)+ \mu  \partial \psi (u_{\mu}(t)) )\Big] \, dt\\
 &&+\ell (u_\mu(0)- u_\mu(T), \frac { u_\mu(T)+ u_\mu(0)}{2})\\
 &=& I_\mu (u_\mu)=0.
\end{eqnarray*}
Since $\Lambda +K$ is regular, $\partial \psi (u_\mu)$ is uniformly bounded and $L$ is weakly lower semi-continuous on $X\times X^*$, we get by letting $\mu\to 0$ that
\begin{eqnarray*}
\ell ( u(T)-u(0), \frac {u(T)+u(0)}{2})+ \int_0^T \Big [\langle \Lambda u(t)+ K u(t), u (t)\rangle + L  (t,u(t),\dot u(t)+\Lambda u(t)+ K u(t))\Big] \, dt \leq 0. 
\end{eqnarray*}
The reverse inequality is true for any $u\in {\cal X}_{p,q}$ since $L$ and $\ell$ are anti-selfdual Lagrangians. \\

\noindent{\bf Proof of Lemma \ref{pert-prop} under condition (A):} Note first that condition (\ref{L.C}) implies that there is a $D>0$ such that 
\begin{equation}\label{L.D}
\hbox{$\int_0^TL(t, u(t), p(t)) dt \geq D (\|p\|^{^s}_{L^q_{X^*}}-1)$ for every $p\in L^q_{X^*}$,}
\end{equation}
where $\frac{1}{r}+\frac{1}{s}=1$. 

However, since $L$ is not supposed to satisfy condition (\ref{L.A}),  we first replace it by its $\lambda$-regularization $L^1_{\lambda}$ which  satisfies all properties of Lemma  \ref{app-lem}.   Therefore,  there exists  $u_{\mu, \lambda } \in {\cal X}_{p,q}$  satisfying
\begin{eqnarray}\label{cdn-A-1}
\dot u_{\mu, \lambda} (t)+\Lambda u_{\mu, \lambda} (t)+ K u_{\mu, \lambda} (t)+ \mu  \partial \psi (u_{\mu, \lambda} (t)))= -\bar \partial   L^1_{\lambda}(t,u_{\mu, \lambda} (t)).
\end{eqnarray}
and
\begin{eqnarray}\label{cdn-A-2}
\ell (u_{\mu, \lambda}(T)-u_{\mu, \lambda}(0),  \frac {u_{\mu, \lambda}(T)+u_{\mu, \lambda}(0)}{2})=\int_0^T \langle \dot u_{\mu, \lambda} (t) ,u_{\mu, \lambda}  (t)\rangle \, dt.
\end{eqnarray}
We shall first  find  bounds for $u_{\mu, \lambda}$ in ${\cal X}_{p,q}$ that are independent of $\mu.$
Multiplying (\ref{cdn-A-1}) by  $u_{\mu, \lambda}$ and integrating, we obtain
\begin{eqnarray}\label{cdn-A-30}
\int _0^T \langle \dot u_{\mu, \lambda} (t)+\Lambda u_{\mu, \lambda} (t)+ K u_{\mu, \lambda} (t)+ \mu  \partial \psi (u_{\mu, \lambda} (t)),u_{\mu, \lambda} (t) \rangle \, dt  = -\int _0^T \langle \bar \partial   L^1_{\lambda}(t,u_{\mu, \lambda} (t)), u_{\mu, \lambda}(t) \rangle \, dt.
\end{eqnarray}
Since $\bar \partial   L^1_{\lambda}(t,.)$ is a {\it maximal monotone } operator, we have 
$\int _0^T \langle \bar \partial   L^1_{\lambda}(t,u_{\mu, \lambda} (t))-\bar \partial   L^1_{\lambda}(t,0), u_{\mu, \lambda} (t)-0 \rangle \, dt \geq 0,$
and therefore
\begin{eqnarray}\label{cdn-A-4}
\int _0^T \langle\bar \partial   L^1_{\lambda}(t,u_{\mu, \lambda} (t)),u_{\mu, \lambda} (t) \rangle \, dt  \geq \int _0^T \langle \bar \partial   L^1_{\lambda}(t,0), u_{\mu, \lambda} (t) \rangle \, dt.
\end{eqnarray}
Taking into account (\ref{cdn-A-2}),  (\ref{cdn-A-4}) and the fact that $\int_0^T\partial \psi (u_{\mu} (t)),u_{\mu} (t) \rangle \geq 0$, it follows from (\ref{cdn-A-30}) that
\begin{eqnarray*}
\ell (u_{\mu, \lambda}(0)-u_{\mu, \lambda}(T),  \frac {u_{\mu, \lambda}(T)+u_{\mu, \lambda}(0)}{2})+\int _0^T \langle \Lambda u_{\mu, \lambda} (t)+ K u_{\mu, \lambda} (t),u_{\mu, \lambda} (t) \rangle \, dt  \leq  -\int _0^T \langle \bar \partial   L^1_{\lambda}(t,0), u_{\mu, \lambda}(t) \rangle \, dt.
\end{eqnarray*}
This implies $\{u_{\mu, \lambda}\}_{\mu}$ is bounded in $L^p_X$, and by the same argument as  under condition (B),  one can prove that $\{\dot u_{\mu, \lambda}\}_{\mu}$ is also bounded in $L^q_{X^*}.$   Consider $u_\lambda \in {\cal X}_{p,q}$ such that  $u_{\mu, \lambda} \rightharpoonup u_{\lambda}$
weakly in ${L^p_{X}}$ and  $\dot {u}_{\mu, \lambda } \rightharpoonup \dot u_{\lambda}$ in ${L^q_{X^*}}$. It follows just like in the proof under condition (B) that 
\begin{eqnarray}
\int_0^T \Big [\langle \Lambda u_{\lambda}(t)+ K u_{\lambda}(t), u_{\lambda}(t) \rangle + L^1_{\lambda} (t,u_{\lambda}(t),\dot u_{\lambda}(t)+
 \Lambda u_{\lambda}(t)+ K u_{\lambda}(t) )\Big] \, dt\nonumber \\+\ell (u_\lambda(0)- u_\lambda(T), \frac { u_\lambda (T)+ u_\lambda (0)}{2})=0,  \label{cdn-A-5}
\end{eqnarray}
and therefore
\begin{eqnarray}\label{cdn-A-55}
\dot u_{\lambda}(t)+
 \Lambda u_{\lambda}(t)+ K u_{\lambda}(t) \in -\bar \partial L^1_{\lambda} (t,u_{\lambda}(t)).
\end{eqnarray}
Now we obtain  estimates on $u_{\lambda}$ in ${\cal X}_{p,q}$.  Since $\ell$ and $L^1_{\lambda}$ are
 bounded from below, it follows from (\ref{cdn-A-5}) that  $\int_0^T \Big [\langle \Lambda u_{\lambda}(t)+ K u_{\lambda}(t), u_{\lambda}(t) \rangle \, dt $  is bounded  and therefore  $u_{\lambda}$ is bounded in $L^p_X$ since 
 \begin{eqnarray*}
 \int_0^T \Big [\langle \Lambda u_{\lambda}(t)+ K u_{\lambda}(t), u_{\lambda}(t) \rangle \, dt
  &\geq& -C(\|u\|_{L^p_X} +1) -\int_0^T \langle \bar \partial L(t, u (t)), u(t)\rangle dt +\int_0^T \langle K u (t), u (t)\rangle \,dt\\
 &\geq&  -C(\|u\|_{L^p_X} +1)-w(\|u\|_{L^p_X})\|u\|_{L^p_X} +  w(\|u\|_{L^p_X}) \|u\|_{L^p_X}^2+ \|u\|_{L^p_X}^{p+1}.
 \end{eqnarray*}
 Setting  $v_{\lambda}(t):=\dot u_{\lambda}(t)+
 \Lambda u_{\lambda}(t)+ K u_{\lambda}(t)$, we get from (\ref{cdn-A-55}) that  
 \[
 -v_{\lambda}(t) = \bar \partial L^1_{\lambda} (t,u_{\lambda}(t))
   = \bar \partial L (t,u_{\lambda}(t)+ \lambda^{q-1} \|v_{\lambda}(t)\|_*^{q-2}J^{-1}v_{\lambda}(t)).
   \] 
   This  together with (\ref{cdn-A-5}) implies that
\begin{eqnarray}
&\int_0^T \Big [\langle \Lambda u_{\lambda}(t)+ K u_{\lambda}(t), u_{\lambda}(t) \rangle +\lambda \|v_{\lambda}(t)\|^q+ L (t,u_{\lambda}(t)+ \lambda \|v_{\lambda}(t)\|_*^{q-2}J^{-1}v_{\lambda}(t),\dot u_{\lambda}(t)+
 \Lambda u_{\lambda}(t)+ K u_{\lambda}(t) )\Big] \, dt \nonumber \\ 
&+ \ell (u_{\lambda}(0)-u_{\lambda}(T),  \frac {u_{\lambda}(T)+
 u_{\lambda}(0)}{2})=0. \label{cdn-A-6}
\end{eqnarray}
It follows  that
$\int_0^T L (t,u_{\lambda}(t)+ \lambda \|v_{\lambda}(t)\|_*^{q-2}J^{-1}v_{\lambda}(t),\dot u_{\lambda}(t)+
 \Lambda u_{\lambda}(t)+ K u_{\lambda}(t) )\, dt$ 
is  bounded from above. 

In view of  (\ref{L.D}), 
 there exists then a
constant $C>0$ such that
\begin{eqnarray}
 \|\dot u_{\lambda}(t)+
 \Lambda u_{\lambda}(t)+ K u_{\lambda}(t)\|_{L^q_{X^*}} \, dt \leq C.
\end{eqnarray}
It follows from the above  that
\begin{eqnarray*}
\|\dot u_{\lambda}\|_{L^q_{X^*}}\leq  \|\Lambda u_{\lambda}\|_{L^q_{X^*}}+
\|K u_{\lambda}\|_{L^q_{X^*}} +C \leq k\|\dot u_{\lambda}\|_{L^q_{X^*}}+ w(\|u\|_{L^p_X})+\|K u_{\lambda}\|_{L^q_{X^*}}
\end{eqnarray*}
from which we obtain
\begin{eqnarray*}
(1-k)\|\dot u_{\lambda}\|_{L^q_{X^*}}\leq  w(\|u_{\lambda}\|_{L^p_X})+\|K u_{\lambda}\|_{L^q_{X^*}}, 
\end{eqnarray*}
which means that $\|\dot u_{\lambda}\|_{L^q_{X^*}}$ is bounded. By letting $\lambda$ go to zero in (\ref{cdn-A-6}), we obtain
\begin{eqnarray*}
\ell ( u(0)-u(T), \frac {u(T)+u(0)}{2})+ \int_0^T \Big [\langle \Lambda u(t)+ K u(t), u (t)\rangle + L  (t,u(t),\dot u(t)+\Lambda u(t)+ K u(t))\Big] \, dt =0
\end{eqnarray*}
where $u$ is a weak limit of $(u_{\lambda})_\lambda $ in ${\cal X}_{p,q}$. 
\hfill $\square$\\

\noindent{\bf Proof of Theorem \ref{main.20}:} First we assume that $\ell$ satisfies condition (\ref{ell.A}), and we shall work towards eliminating the perturbation $K$. Let $L^2_{\lambda}$ be the $\lambda-$regularization of $L$ with respect to the second variable, in such a way that $L^2_{\lambda}$
satisfies (\ref{coercive-2}). Indeed
\begin{eqnarray}\label{560}
\int_0^T \langle \bar \partial L^2_{\lambda} (t,u(t)) + \Lambda u(t), u(t) \rangle \, dt &=& \int_0^T \langle \bar \partial L(t, u(t)) + \Lambda u(t)
 + \lambda^{p-1}\|u\|^{p-2} Ju(t), u(t) \rangle \, dt \nonumber  \\ & \geq & \int_0^T \langle \bar \partial L (t,u(t)) + \Lambda u(t) , u(t) \rangle \, dt \geq -C \|u\|_{L^p_X}.
\end{eqnarray}
Moreover, we have in view of (\ref{L.C}) that 
\begin{eqnarray}\label{L.E}
\int_0^TL^2_\lambda (t, u,p)\, dt \geq  -D+\frac{\lambda^{p-1}}{p}\|u\|_{L^p_X}^p. 
\end{eqnarray}

From Lemma \ref{pert-prop}, we get  for each $\epsilon >0$,  $u_{\epsilon, \lambda} \in {\cal X}_{p,q}$ such that
\begin{eqnarray}\label{main-2-1}
 \int_0^T \Big [\langle \Lambda u_{\epsilon, \lambda}(t)+\epsilon  K
u_{\epsilon, \lambda}(t), u_{\epsilon, \lambda}(t) \rangle + L^2_{\lambda}  (t,u_{\epsilon, \lambda}(t),\dot u_{\epsilon, \lambda}(t)+\Lambda u_{\epsilon, \lambda}(t)+ \epsilon K u_{\epsilon, \lambda}(t)) \Big ] \, dt\
\nonumber\\
+\ell (u_{\epsilon, \lambda} (0)- u_{\epsilon, \lambda} (T), \frac {u_{\epsilon, \lambda} (T)+ u_{\epsilon, \lambda} (0)}{2})=0,  
\end{eqnarray}
and 
\begin{equation}\label{new}
 \dot u_{\epsilon, \lambda} (t)+\Lambda u_{\epsilon, \lambda} (t)+\epsilon  K u_{\epsilon, \lambda} (t) \in  - \bar \partial   L^2_{\lambda}(t,u_{\epsilon, \lambda} (t)).
\end{equation}
We shall first  find  bounds for $u_{\epsilon, \lambda}$ in ${\cal X}_{p,q}$ that are independent of $\epsilon.$
Multiplying (\ref{new}) by  $u_{\epsilon, \lambda}$ and integrating, we obtain
\begin{eqnarray}\label{cdn-A-3}
\int _0^T \langle \dot u_{\epsilon, \lambda} (t)+\Lambda u_{\epsilon, \lambda} (t)+\epsilon  K u_{\epsilon, \lambda} (t),u_{\epsilon, \lambda} (t) \rangle \, dt  = -\int _0^T \langle \bar \partial   L^2_{\lambda}(t,u_{\epsilon, \lambda} (t)), u_{\epsilon, \lambda}(t) \rangle \, dt.
\end{eqnarray}
It follows from (\ref{560}) and the above equality that
\begin{eqnarray}\label{cdn-B-30}
 \int_0^T
\langle \dot u_{\epsilon, \lambda} (t)+ \epsilon K u_{\epsilon, \lambda}  (t),u_{\epsilon, \lambda}  (t) \rangle  \leq C\|u_{\epsilon, \lambda}\|_{L^p_X}, 
\end{eqnarray}
and therefore
\begin{eqnarray*}
\ell (u_{\epsilon, \lambda} (0)-u_{\epsilon, \lambda} (T),  \frac {u_{\epsilon, \lambda} (T)+u_{\epsilon, \lambda} (0)}{2})+ \int_0^T
\langle  \epsilon K u_{\epsilon, \lambda}   (t),u_{\epsilon, \lambda}  (t) \rangle  \leq C\|u_{\epsilon, \lambda} \|_{L^p_X}, 
\end{eqnarray*}
which in view of (\ref{main-2-1}) implies that 
\begin{eqnarray*}
| \int_0^T  L^2_{\lambda}  (t,u_{\epsilon, \lambda}(t),\dot u_{\epsilon, \lambda}(t)+\Lambda u_{\epsilon, \lambda}(t)+ \epsilon K u_{\epsilon, \lambda}(t))  \, dt| 
&\leq&C\|u_{\epsilon, \lambda} \|_{L^p_X}.
\end{eqnarray*}
By (\ref{L.E}), we deduce that   
$\{u_{\epsilon, \lambda}\}_{\mu}$ is bounded in $L^p_X$.  The same reasoning as above then shows   that $\{\dot u_{\epsilon, \lambda}\}_{\mu}$ is also bounded in $L^q_{X^*}$.  
Again, the regularity of  $\Lambda$ and the lower semi-continuity of $L$,  yields the existence of $u_{\lambda} \in {\cal X}_{p,q}$ such that
\begin{eqnarray}\label{main-2-2}
\ell ( u_{\lambda}(0)-u_{\lambda}(T), \frac {u_{\lambda}(T)+u_{\lambda}(0)}{2})+ \int_0^T \Big [\langle \Lambda u_{\lambda}(t)
, u_{\lambda}(t) \rangle + L^2_{\lambda}  (t,u_{\lambda}(t),\dot u_{\lambda}(t)+\Lambda u_{\lambda}(t)) \Big ] \, dt =0.
\end{eqnarray}
In other words, 
\begin{eqnarray}
& \int_0^T \Big [\langle \Lambda u_{\lambda}(t)
, u_{\lambda}(t) \rangle + L  (t,u_{\lambda}(t),\dot u_{\lambda}(t)+\lambda^{p-1} \|u_{\lambda}(t)\|^{p-2} u_{\lambda}(t)+\lambda J u_{\lambda}(t)) +\lambda^{p-1} \|u_{\lambda}(t)\|^p\Big ] \, dt\nonumber \\
&+\ell (u_\lambda (0)- u_\lambda (T), \frac { u_\lambda (T)+ u_\lambda (0)}{2})=0. \label{main-2-3}
\end{eqnarray}
Now since $L$ is $\Lambda-$coercive we get that $(u_{\lambda})_\lambda$ is bounded in ${\cal X}_{p,q}.$ Suppose  $u_{\lambda}\rightharpoonup \bar u$ in $L^p_X$
and  $\dot u_{\lambda}\rightharpoonup \dot {\bar u}$ in $L^q_{X^*}.$ It follows from (\ref{Lambda-1}) that $\Lambda u_{\lambda}$
is bounded in $L^q_{X^*}.$
Again, we deduce that
\begin{eqnarray*}
\ell ( \bar u(T)-\bar u(0), \frac {\bar u(T)+\bar u(0)}{2})+ \int_0^T
\Big [\langle \Lambda \bar u(t) , \bar u(t)\rangle + L  (t, \bar u(t),\dot {\bar
u}(t)+\Lambda \bar u(t)) \Big ] \, dt =0.
\end{eqnarray*}
Now, we show that we can do without assuming that $\ell$ satisfies (\ref{ell.A}), but that it is bounded below while $(0,0)\in {\rm Dom}(\ell)$. Indeed, let $\ell_{\lambda}:=\ell^{1,2}_\lambda$ be the
$\lambda$-regularization of the anti-selfdual Lagrangian $\ell$ in both variables. 
Then  $\ell_{\lambda}$ satisfies (\ref{ell.A})
and therefore there exists $x_{\lambda} \in
{\cal X}_{p,q}$ such that
\begin{eqnarray}\label{main-2-4}
\ell_{\lambda}( x_{\lambda}(T)-x_{\lambda}(0), \frac {x_{\lambda}(T)+x_{\lambda}(0)}{2})+ \int_0^T \Big [\langle \Lambda x_{\lambda}(t)
, x_{\lambda} (t)\rangle + L  (t,x_{\lambda}(t),\dot x_{\lambda}(t)+\Lambda x_{\lambda}(t)) \Big ] \, dt =0.
\end{eqnarray}
Since  $\ell$ is bounded from below so is $\ell_{\lambda}$.   This together
with (\ref{main-2-4}) imply that the family $ \int_0^T \Big [\langle \Lambda x_{\lambda}(t)
, x_{\lambda}(t) \rangle + L  (t,x_{\lambda}(t),\dot
x_{\lambda}(t)+\Lambda x_{\lambda}(t)) \Big ] \, dt$ is bounded
above. Again, since $L$ is $\Lambda-$coercive, we obtain that $(x_{\lambda})_\lambda$
is bounded in ${\cal X}_{p,q}$.  The continuity of the injection  ${\cal X}_{p,q}
\subseteq C([0,T];H)$ also ensures  the boundedness
 of
$(x_{\lambda}(T))_\lambda$ and $(x_{\lambda}(0))_\lambda$ in $H$.  Consider 
$\bar x \in {\cal X}_{p,q}$ such that  
$x_{\lambda}\rightharpoonup \bar x$ in $L^p_X$ and  $\dot
x_{\lambda}\rightharpoonup \dot {\bar x}$ in $L^q_{X^*}.$  It
follows from the regularity of $\Lambda$ and the lower semi-continuity of $\ell$ and $L$ 
that
\begin{eqnarray*}
\ell ( \bar x(T)-\bar x(0), \frac {\bar x(T)+\bar x(0)}{2})+ \int_0^T
\Big [\langle \Lambda \bar x(t) , \bar x(t)\rangle + L  (\bar x(t),\dot {\bar
x}(t)+\Lambda \bar x(t)) \Big ] \, dt =0,
\end{eqnarray*}
and therefore $\bar x$ satisfies equation (\ref{principle3}).

\section{Application to Navier-Stokes evolutions}

The most basic time-dependent  selfdual Lagrangians are of the
form
$
L(t,x,p)=\varphi (t, x) +\varphi^{*}(t, -p)
$
where for each $t$, the function $x\to \varphi (t, x)$ is convex and lower semi-continuous on $X$. Let  now $\psi: H \rightarrow \R \cup\{+\infty\}$
be another convex lower semi-continuous function which is bounded from below and such that $0 \in {\rm Dom} (\psi)$, and set $\ell (a,b)= \psi (a)+ \psi^* (-b)$.
The above principle then yields that  if  for some $C_1, C_2>0$, we have 
\[
\hbox{$C_1\big(\|
x\|_{L^p_X}^p-1\big)\leq\int_0^T\phi\big( t,x(t)\big)\,dt\leq C_2\big(\|
x\|_{L^p_X}^p+1\big)$, for all $x\in L^p_X$}, 
\]
then for every  regular map $\Lambda$  satisfying (\ref{Lambda-1}) and either one of conditions  (A) or (B) 
in Theorem \ref{main.20},  the infimum of the functional
\[
I(x)=\int_0^T \big [\phi( t, x(t))+\phi^*(t, -\dot{x}(t)-\Lambda x(t))+ \langle \Lambda  x(t)
,  x(t) \rangle \big ] \,dt + \psi (x(0)-x(T)) +\psi^*(- \frac{x(0)+ x(T)}{2})
\]
on ${\cal X}_{p,q}$ is zero and is attained at a solution $x(t)$ of the following equation
 \begin{eqnarray*}
 \left \{ \begin{array}{lcl}
\hfill -\dot x(t)-\Lambda x(t) &\in & \partial\phi\big( t,x(t)\big)  \quad {\rm for\,  all}\, t\in [0,T] \\
 \hfill -\frac{ x(0)+ x(T)}{2} & \in & \partial \psi (x(0)-x(T)).
\end{array}\right.
\end{eqnarray*}
As noted in the introduction, the boundary condition above 
is quite general and it includes as particular case the more traditional ones such as initial-value problems, periodic and anti-periodic orbits.  It suffices  to choose $\ell (a,b)= \psi (a)+ \psi^* (-b)$ accordingly.
  \begin{itemize}
\item For the initial boundary condition $x(0)=x_0$ for a given $x_0\in
H$, we choose $\psi (x)= \frac{1}{4}\|x\|_H^2-\braket{x}{x_0}$.
\item For periodic solutions $x(0)=x(T)$,  $\psi$ is chosen as:
\begin{eqnarray*}\psi(x)=\left\{\begin{array}{ll}
0 \quad &x=0\\
+\infty &\mbox{elsewhere}.\end{array}\right.
\end{eqnarray*}
\item For anti-periodic solutions $x(0)=-x(T)$, it suffices to choose $\psi (x)=0$ for each $x \in H.$
\end{itemize}

As a consequence of the above theorem, we provide a variational resolution to evolution equations involving  nonlinear operators such as the Navier-Stokes equation with various boundary conditions:
  \begin{equation}
\label{T-NS}
 \left\{ \begin{array}{lcl}
    \hfill
 \frac {\partial u}{\partial t}+(u\cdot \nabla)u +f &=&\nu \Delta u - \nabla  p \quad \hbox{\rm on $ \Omega$},\\
\hfill {\rm div} \, u&=&0 \quad  \quad  \quad  \quad  \quad \hbox{\rm on  $[0, T]\times \Omega$},\\
\hfill u&=&0 \quad  \quad  \quad  \quad  \quad \hbox{\rm on $[0, T]\times \partial \Omega$}, 
\end{array}\right.
\end{equation}
where $\Omega$ is a smooth domain of $\R^n$,  $f\in L^{2}_{X^*}([0,T])$, $\nu>0$.

 Indeed,  setting  $X=\{u\in H^{1}_0(\Omega; {\bf R}^{n}); {\rm div} v=0\}$, and $H=L^2(\Omega)$, 
  we write the above problem in the form
 \begin{eqnarray}
\label{T-NS-eq}
 \left\{ \begin{array}{lcl}
  \hfill \frac{\partial u}{\partial t}+ \Lambda u   &\in& -\partial \Phi (t,u) \\
 \hfill  \frac{ u(0)+ u(T)}{2} & \in& -\bar \partial \ell (u(0)-u(T)), 
\end{array}\right.
\end{eqnarray}
where $\ell$ is any anti-self dual Lagrangian on $H \times H$, while the convex functional  $\Phi$  and the nonlinear operator $\Lambda$ are defined by: 
 \begin{equation}
\hbox{$\Phi (t,u)=\frac{\nu}{2}
\int_{\Omega}\Sigma_{j,k=1}^{3}(\frac {\partial u_{j} } {\partial
x_{k}})^{2}\, dx + \langle u,f (t,x) \rangle $ and   
$\Lambda u:= (u\cdot
\nabla)u$.}
\end{equation}
 Note that $\Lambda: X\to X^*$ as long as the dimenison $N\leq 4$. On the other hand, when $\Lambda$ lifts to path space, we have the following 
\begin{lemma}(1) When  $N=2 $,  the operator $\Lambda : {\cal
X}_{2,2}\rightarrow L^2_{X^*}$ is regular.\\
(2) When  $N=3$,  the operator $\Lambda$ is regular from  ${\cal X}_{4,\frac{4}{3}}\rightarrow
L^{\frac{4}{3}}_{X^*}$  as well as from $ {\cal X}_{2,\frac{4}{3}}\cap
L^{\infty}(0,T;H)$ to  $L^{\frac{4}{3}}_{X^*}$.
\end{lemma}
\noindent{\bf Proof:} First note that the three embeddings ${\cal X}_{2,2}
\subseteq L^2_H,$ $ {\cal X}_{4,\frac{4}{3}}\subseteq L^2_H,$ and $ {\cal
X}_{2,\frac{4}{3}}\subseteq L^2_H,$ are compact.\\
Assume that $ N=3$, let  $u^{n}\to u$ weakly in ${\cal X}_{4,\frac{4}{3}}$,
and fix $v\in C^1([0,T] \times \Omega).$  We have that
\[
\int_0^T \langle \Lambda u^{n}, v\rangle=\int_0^T
\int_{\Omega}\Sigma_{j,k=1}^{3}u^{n}_{k}\frac {\partial u_{j}^n } {\partial
x_{k}}v_{j}\, dx \, dt=-\int_0^T \int_{\Omega}\Sigma_{j,k=1}^{3}u^{n}_{k}\frac {\partial v_{j} } {\partial x_{k}}u^n_{j}\, dx.
\]
Therefore 
\begin{eqnarray}\label{sequence}
\big |\int_0^T \langle \Lambda u^{n}-\Lambda u, v\rangle \big | & =&\big
|\Sigma_{j,k=1}^{3}\int_0^T \int_{\Omega}(u^{n}_{k}\frac {\partial v_{j} }
{\partial x_{k}}u^n_{j}-u_{k}\frac {\partial v_{j} } {\partial
x_{k}}u_{j})\, dx \, dt \big| \nonumber \\
& \leq & \|v\|_{C^1([0,T] \times \Omega)} \Sigma_{j,k=1}^{3}\int_0^T
\int_{\Omega}\big  |u^{n}_{k}u^n_{j}-u_{k}u_{j}\big | \, dx \, dt.
\end{eqnarray}
Also
\begin{eqnarray}\label{sequence1}
\int_0^T \int_{\Omega}|u^{n}_{k}u^n_{j}-u_{k}u_{j}| \, dx \, dt & \leq &
\int_0^T \int_{\Omega}|u^{n}_{k}u^n_{j}-u_{k}u^n_{j}| \, dx \, dt +\int_0^T
\int_{\Omega}|u_{k}u^n_{j}-u_{k}u_{j}| \, dx \, dt \nonumber \\   & \leq &
\|u^n_{j}\|_{L^2_H} \|u^{n}_{k}-u_{k}\|_{L^2_H} + \|u_k\|_{L^2_H}
\|u^{n}_{j}-u_{j}\|_{L^2_H} \rightarrow 0.
\end{eqnarray}
 Moreover, we have for $N=3$, the following standard estimate (\cite{Te})
\begin{eqnarray}\label{estimate}
\|\Lambda u^n\|_{X^*}\leq c |u^n|_H^{\frac{1}{2}}\|u^n\|_X^{\frac{3}{2}}.
\end{eqnarray}
Since ${\cal X}_{4,\frac{4}{3}}\subseteq C(0,T;H)$ is continuous, we
obtain
\begin{eqnarray}\label{estimate1}
\|\Lambda u^n\|_{L^{\frac{4}{3}}_{X^*}}\leq c |u^n|^{\frac{1}{2}}_{C(0,T;H)}\|u^n\|_{L^2_X}^{\frac{3}{4}}\leq c \|u^n\|^{\frac{1}{2}}_{{\cal X}_{4,\frac{4}{3}}}\|u^n\|_{L^2_X}^{\frac{3}{4}}
\end{eqnarray}
from which we conclude that $\Lambda u^n$ is a bounded sequence in
$L^{\frac{4}{3}}_{X^*}$, and therefore  the convergence holds for each $v \in
L^4_X.$

Now, since ${\cal X}_{2,2}\subseteq C(0,T;H)$ is also continuous, the same argument
works for $N=2$, the only difference being that we have the following estimate which is better that (\ref{estimate}), 
 \begin{eqnarray}\label{nick.1}
\|\Lambda u^n\|_{X^*}\leq c |u^n|_H\|u^n\|_X.
\end{eqnarray}

To consider the case $\Lambda: {\cal X}_{2,\frac{4}{3}}\cap L^{\infty}(0,T;H) \to L^{\frac{4}{3}}_{X^*}$, we note that 
  relations (\ref{sequence}) and
(\ref{sequence1}) still hold if $u_n\to u$ weakly in ${\cal X}_{2,\frac{4}{3}}$. We also have estimate  (\ref{estimate}). However, unlike the above,  one cannot deduce (\ref{estimate1})  since we do not have necessarily a continuous embbeding from ${\cal X}_{2,\frac{4}{3}}\subseteq C(0,T;H)$. However, if $(u_n)$ is also assumed to be bounded in $L^{\infty}(0,T;H)$, then  we get the following estimate from (\ref{estimate}),
\begin{eqnarray}\label{nick.2}
\|\Lambda u^n\|_{L^{\frac{4}{3}}_{X^*}}\leq c |u^n|^{\frac{1}{2}}_{L^{\infty}(0,T;H)}\|u^n\|_{L^2_X}^{\frac{3}{4}}
\end{eqnarray}
which ensures the boundedness of $\Lambda u^n$ in $L^{\frac{4}{3}}_{X^*}.$
\hfill $\square$ 
 
We now prove Corollaries \ref{T-NS2} and \ref{T-NS3} stated in the introduction.\\

  \noindent{\bf Proof of Corollary \ref{T-NS2}:}   By the preceeding lemma, one can verify that the operator  $\Lambda : {\cal X}_{2,2} \rightarrow L^2_{X^*}$
  satisfies condition  (\ref{Lambda-10})  and (\ref{Lambda-20}).
 Therefore  the infimum of the functional
\[
I(u)=\int_0^T \big [\Phi( t, u(t))+\Phi^*(t, -\dot{u}(t)- (u\cdot \nabla)u(t)) \big ] \,dt +\ell \big ( u(0)-u(T), \frac{u(0)+ u(T)}{2}  \big )
\]
on ${\cal X}_{2,2}$ is zero and is attained at a solution $u(t)$ of (\ref{T-NS}). $\square$\\

  \noindent{\bf Proof of Corollary \ref{T-NS3}:}   We start by considering  the following functional on the space ${\cal X}_{4, \frac {4}{3}}$. 
\[
I_{\epsilon}(u):=\int_0^T \big [\Phi_{\epsilon}( t, u(t))+\Phi_{\epsilon}^*(t, -\dot{u}(t)- (u\cdot \nabla)u(t)) \big ] \,dt + \ell (u(0)-u(T), \frac{u(0)+ u(T)}{2})
\]
where $\Phi_{\epsilon}( t, u)= \Phi( t, u)+ \frac {\epsilon}{4} \|u\|^4$. In view of the preceeding lemma, the operator  $\Lambda u:=(u\cdot \nabla)u$ and $\Phi_{\epsilon}$
satisfy all properties of Theorem \ref{main.20}. In particular, we have the estimate
\begin{equation}\label{NS.estimate}
 \|\Lambda u\|_{X^*} \leq c |u|_H^{1/2} \|u\|_X^{3/2} \quad \text { for every } u \in X. 
 \end{equation} 
 It follows from Theorem \ref{main.20}, that there exists  $u_{\epsilon } \in {\cal X}_{4, \frac {4}{3}} $
with $I_{\epsilon}(u_{\epsilon})=0$. This implies that 
\begin{eqnarray}\label{T-NS3-1}
 \left\{ \begin{array}{lcl}
    \hfill
\frac{\partial u_{\epsilon}}{\partial t}+(u_{\epsilon}\cdot \nabla)u_{\epsilon} +f(t,x) &=&\nu \Delta u_{\epsilon}+ {\epsilon} \|u_{\epsilon}\|^2 \Delta u_{\epsilon}- \nabla  p_{\epsilon} \quad \hbox{\rm on $[0,T]\times \Omega$}\\
\hfill {\rm div} u_{\epsilon}&=&0 \quad \quad \quad \hbox{\rm on $[0,T]\times \Omega$}\\
\hfill u_\epsilon &=&0 \quad \quad \quad \hbox{\rm on $[0, T]\times \partial \Omega$}.\\
\hfill -\frac{ u_{\epsilon}(0)+ u_{\epsilon}(T)}{2} & =& \bar \partial \ell (u_{\epsilon}(0)-u_{\epsilon}(T)).
\end{array}\right.
\end{eqnarray}

Now, we show that $(u_{\epsilon})_\epsilon$ is bounded in ${\cal X}_{2, 4/3}.$
Indeed, multiply (\ref{T-NS3-1}) by $u_{\epsilon}$  to get
\begin{eqnarray*}
\frac {d}{dt} \frac {|u_{\epsilon}(t)|^2} {2}+\nu \|u_{\epsilon} (t)\|^2+ \epsilon \|u_{\epsilon} (t)\|^4  =  \langle f(t), u_{\epsilon} (t) \rangle
 \leq  \frac {\nu}{2} \|u_{\epsilon} (t)\|^2 +\frac {2} {\nu} \|f(t)\|_{X^*}^2
\end{eqnarray*}
so that
\begin{eqnarray}\label{T-NS3-2}
\frac {d}{dt} \frac {|u_{\epsilon}(t)|^2} {2}+ \frac {\nu}{2} \|u_{\epsilon} (t)\|^2+ \epsilon \|u_{\epsilon} (t)\|^4
 \leq  \frac {2} {\nu} \|f(t)\|_{X^*}^2.
\end{eqnarray}
Integrating (\ref{T-NS3-2}) over $[0,s],$  $(s <T)$ we obtain
\begin{eqnarray}\label{T-NS3-3}
 \frac {|u_{\epsilon}(s)|^2} {2}- \frac {|u_{\epsilon}(0)|^2} {2}+ \frac {\nu}{2} \int _0^s \|u_{\epsilon} (t)\|^2+ \epsilon \int _0^s  \|u_{\epsilon} (t)\|^4
 \leq  \frac {2} {\nu}  \int _0^s  \|f(t)\|_{X^*}^2.
\end{eqnarray}

On the other hand, it follows from (\ref{T-NS3-1}) that
$\ell (u_{\epsilon}(0)-u_{\epsilon}(T), \frac{ u_{\epsilon}(0)+ u_{\epsilon}(T)}{2}).=\frac {|u_{\epsilon}(T)|^2} {2}- \frac {|u_{\epsilon}(0)|^2} {2}.$
 Considering this together with (\ref{T-NS3-3}) with $s=T$, we get
\begin{eqnarray}\label{T-NS3-4}
\ell (u_{\epsilon}(0)-u_{\epsilon}(T), \frac{ u_{\epsilon}(0)+ u_{\epsilon}(T)}{2}) + \frac {\nu}{2} \int _0^T \|u_{\epsilon} (t)\|^2+ \epsilon \int _0^T  \|u_{\epsilon} (t)\|^4
 \leq  \frac {2} {\nu}  \int _0^T  \|f(t)\|_{X^*}^2.
\end{eqnarray}

Since $\ell$ is bounded from below and is coercive in both variables,  it follows from the above that $(u_{\epsilon})_\epsilon$ is bounded in
 $L^2_X$, that $(u_{\epsilon}(T))_\epsilon$ and $(u_{\epsilon}(0))_\epsilon$ are bounded in $H$, and that    $\epsilon \int _0^T  \|u_{\epsilon} (t)\|^4$  is also bounded. It also 
  follows from  (\ref{T-NS3-3}) coupled with the boundedness of $(u_{\epsilon}(0))_\epsilon$, that $u_{\epsilon}$ is bounded in $L^{\infty}(0,T; H)$.   
 Estimate (\ref{NS.estimate}) combined with the boundedness of $(u_{\epsilon})_\epsilon$ in $L^{\infty}(0,T; H) \cap L^2_X$ implies that $(\Lambda u_{\epsilon})_\epsilon$ is
bounded in $L^{4/3}_X$. We also have the estimate
\[\big \|\nu \Delta u_{\epsilon}+ {\epsilon} \|u_{\epsilon}\|^2 \Delta u_{\epsilon} \big \|_{X^*} \leq \nu \|u_{\epsilon}\| + \epsilon \|u_{\epsilon}\|^3\]
 which implies that $\nu \Delta u_{\epsilon}+ {\epsilon} \|u_{\epsilon}\|^2 \Delta u_{\epsilon}$ is bounded in $L^{4/3}_{X^*}.$

It also follows from (\ref{T-NS3-1}) that for each $v \in L^4_X$,  we have
\begin{eqnarray}\label{T-NS3-5}
\int_0^T \langle \frac{\partial u_{\epsilon}}{\partial t}, v \rangle \, dt = \int_0^T \langle -(u_{\epsilon}\cdot \nabla)u_{\epsilon} -f(t,x)
+\nu \Delta u_{\epsilon}+ {\epsilon} \|u_{\epsilon}\|^2 \Delta u_{\epsilon}, v \rangle \, dt.
\end{eqnarray}
Since the right hand side is uniformly bounded with respect to $\epsilon$, so is the left hand side, which implies that 
 $ \frac{\partial u_{\epsilon}}{\partial t}$ is bounded in $L^{4/3}_{X^*}.$ Therefore, there exists $u \in {\cal X}_{2, 4/3}$ such that
\begin{eqnarray}
\hfill u_{\epsilon} & \rightharpoonup & u \text{   weakly in   } L^2_X,\label{T-NS3-6}\\
\hfill  \frac{\partial u_{\epsilon}}{\partial t}   & \rightharpoonup  & \frac {\partial u_{\epsilon}}{\partial t}   \text  { weakly in } L^{4/3}_{X^*},\label{T-NS3-7}\\
 {\epsilon} \|u_{\epsilon}\|^2 \Delta u_{\epsilon}  & \rightharpoonup  &  0  \text  { weakly in } L^{4/3}_{X^*},\label{T-NS3-8}\\
\hfill u_{\epsilon} (0)  & \rightharpoonup & u (0)   \text { weakly in } H,\label{T-NS3-9}\\
\hfill  u_{\epsilon} (T)  & \rightharpoonup & u(T)   \text { weakly in } H. \label{T-NS3-10}
\end{eqnarray}
Letting $\epsilon $ approach to zero in  (\ref{T-NS3-5}), it follows from (\ref{T-NS3-6})-(\ref{T-NS3-10})
that
\begin{eqnarray}\label{T-NS3-11}
\int_0^T \langle \frac{\partial u}{\partial t}, v \rangle \, dt = \int_0^T \langle -(u\cdot \nabla)u-f(t,x).
+\nu \Delta u, v \rangle \, dt.
\end{eqnarray}
Also it follows from (\ref{T-NS3-9}), (\ref{T-NS3-10}) and   (\ref{T-NS3-1})  and the fact that $\bar \partial \ell$ is maximal monotone that
\begin{eqnarray}\label{T-NS3-12}
-\frac {u(0)+u(T)}{2}=\bar \partial \ell (u(0)-u(T)).
\end{eqnarray}
  (\ref{T-NS3-11}) and (\ref{T-NS3-12}) yield that $u$ is a weak solution of
\begin{eqnarray}\label{T-NS3-13}
 \left\{ \begin{array}{lcl}
    \hfill
\frac{\partial u}{\partial t}+(u\cdot \nabla)u +f(t,x) &=&\nu \Delta u - \nabla  p \quad \hbox{\rm on $[0,T]\times \Omega$},\\
\hfill {\rm div} u&=&0 \quad \quad \quad \quad \quad \hbox{\rm on $[0,T]\times \Omega$},\\
\hfill u&=&0 \quad \quad \quad \quad \quad  \hbox{\rm on $[0, T]\times \partial \Omega$}.\\
\hfill -\frac {u(0)+u(T)}{2} & =& \bar \partial \ell (u(0)-u(T)) \quad \hbox{\rm on $ \Omega$.}
\end{array}\right.
\end{eqnarray}
Now we prove inequality (\ref{T-NS3-in}). Since $I_{\epsilon}(u_{\epsilon})=0$,  a standard argument
 (see the proof of Theorem \ref{main.20}) yields
 that
$I(u)\leq \liminf_{\epsilon}I_{\epsilon}(u_{\epsilon})=0$,
thereby giving that 
\[
I_{\epsilon}(u):=\int_0^T \big [\Phi ( t, u(t))+\Phi^*(t, -\dot{u}(t)- (u\cdot \nabla)u(t)) \big ] \,dt + \ell (u(0)-u(T), \frac{u(0)+ u(T)}{2}) \leq 0.
\]
On the other hand it follows from (\ref{T-NS3-12}) that $\ell (u(0)-u(T),- \frac{u(0)+ u(T)}{2}) = \frac {|u(T)|^2}{2}- \frac {|u(0)|^2}{2}$. This
together with the above inequality gives
\[
\frac {|u(T)|^2}{2}+ \int_0^T \big [\Phi ( t, u(t))+\Phi^*(t, -\dot{u}(t)- (u\cdot \nabla)u(t)) \big ] \,dt  \leq \frac {|u(0)|^2}{2}.
\]
\begin{corollary} In dimension  $N=3$, there exists for any given $\alpha$ with $|\alpha |<1$, a weak solution of the equation
 solutions:
\begin{eqnarray*}
 \left\{ \begin{array}{lcl}
    \hfill
\frac{\partial u}{\partial t}+(u\cdot \nabla)u +f(t,x) &=&\nu \Delta u - \nabla  p \quad \hbox{\rm on $[0,T]\times \Omega$},\\
\hfill {\rm div} u&=&0 \quad \quad \quad \hbox{\rm on $[0,T]\times \Omega$},\\
\hfill u&=&0 \quad \quad \quad \hbox{\rm on $[0, T]\times \partial \Omega$}.\\
\hfill u(0) & = &  \alpha u(T).
\end{array}\right.
\end{eqnarray*}
\end{corollary}
\noindent{\bf Proof:}  For each $\alpha $ with  $|\alpha |<1$ there exists $\lambda >0$ such that $\alpha=\frac{\lambda-1}{\lambda+1}.$
 Now  consider $\ell(a,b)=\psi_{\lambda}(a)+\psi^*_{\lambda}(-b)$ where $\psi_{\lambda} (a)= \frac {\lambda}{4} |a|^2.$ $\square$\\
  
  {\bf Navier-Stokes evolutions driven by their boundary:} 
   We now consider the following evolution equation.
 \begin{eqnarray}
 \label{NSE50}
  \left\{ \begin{array}{lcl}
     \hfill
 \frac{\partial u}{\partial t}+(u\cdot \nabla)u +f &=&\nu \Delta u - \nabla  p \quad \hbox{\rm on $[0,T]\times \Omega$}\\
 \hfill {\rm div} u&=&0 \quad \quad \quad  \quad  \quad \hbox{\rm on $[0,T]\times \Omega$}\\
 \hfill u(t,x)&=&u^0(x) \quad  \quad  \quad \hbox{\rm on $[0,T]\times \partial \Omega$}\\
 \hfill u(0,x)&=&\al u(T,x) \quad  \quad \hbox{\rm on $\Omega$} 
 \end{array}\right.
 \end{eqnarray}
   where $\int_{\partial \Omega} u^{0} {\bf \cdot n}\, d\sigma =0$, $\nu>0$ and $f\in L^p_{X^*}$.  Assuming  that $u^{0} \in H^{3/2}(\partial \Omega)$ and that $\partial \Omega$ is connected,  Hopf's extension theorem again yields  the existence of $v^0\in H^{2}(\Omega)$ such that
  \begin{equation}
  \label{Hopf}
 v^0=u^{0}\,\, \hbox{\rm on $\partial \Omega$,\quad  ${\rm div}\, v^0=0$\quad  and \quad  $\int_{\Omega}\Sigma_{j,k=1}^{n}u_{k}\frac {\partial v^0_{j} } {\partial x_{k}}u_{j}\, dx\leq \epsilon \|u\|^{2}_{X}$ for all $u\in X$}
  \end{equation}
  where  $V=\{u\in H^{1}(\Omega; {\bf R}^{n}); {\rm div} u=0\}$.
  Setting $v=u+v^0$, then solving  (\ref{NSE50}) reduces to finding a solution in the path space ${\cal X}_{2,2}$ corresponding to the Banach space $X=\{u\in H^{1}_0(\Omega; {\bf R}^{n}); {\rm div} v=0\}$ and the Hilbert space $H=L^2(\Omega)$ 
   for
   \begin{eqnarray}
  \label{NSE51}
  \frac{\partial u}{\partial t}+ (u\cdot \nabla)u +(v^0\cdot \nabla)u +(u\cdot \nabla)v^0    &\in& -\partial \Phi (u) \\
  \hfill  u(0)- \al u(T)&=& (\al-1)v^0. 
   \nonumber 
  \end{eqnarray}
  where $\Phi (t,u)=\frac{\nu}{2} \int_{\Omega}\Sigma_{j,k=1}^{3}(\frac {\partial u_{j} } {\partial x_{k}})^{2}\, dx+ \langle g, u \rangle,$  and where 
 \[
g:= f -\nu \Delta v^0 + (v^0\cdot \nabla)v^0\in L^p_{V^*}.
\]
In other words, this is an equation of the form
\begin{equation}
 \label{Lax.Mil.again}
\frac{\partial u}{\partial t}+  \Lambda u\in -\partial \Phi (t,u) 
  \end{equation}
  where $\Lambda u:= (u\cdot \nabla)u+(v^0\cdot \nabla)u +(u\cdot \nabla)v^0$ is the nonlinear regular operator  $N=2$ or $N=3.$

Now  recalling the fact  that the component $Bu:=(v^0\cdot \nabla)u$  is skew-symmetric, it follows from   Hopf's estimate that  
  \[
C\|u\|_V^2\geq \Phi (t,u)+\langle \Lambda u,u\rangle \geq (\nu-\epsilon)\|u\|^{2}+\langle g,u  \rangle \quad {\rm for\,  all} \,\, u\in X.
  \]
As in Corollary \ref{T-NS3} we have the following.
\begin{corollary} Assume $N=3.$ Consider  $\ell$ to be a selfdual Lagrangian on $H \times H$ that is coercive in both variables.  Then,
there  exists $u \in {\cal X}_{2, {\frac {4}{3}}} $ such that
\[
I(u)=\int_0^T \big [\Phi( t, u(t))+\Phi^*(t, -\dot{u}(t)- \Lambda u(t)) +\langle u(t), \Lambda u(t) \rangle  \big ] \,dt + \ell (u(0)-u(T), \frac{u(0)+ u(T)}{2})\leq 0
\]
 and  $u$ is  a weak solution of  (\ref{NSE50}).  
\end{corollary}
To obtain boundary condition given in (\ref{NSE51}) that is $u(0)-\al u(T)= (\al-1)v^0,$ consider $\ell(a,b)= \psi_{\lambda}(a)+\psi^*_{\lambda}(-b)$ where $\al= \frac{\lambda-1}{\lambda+1}$ and $\psi_{\lambda}(a)=\frac {\lambda}{4}|a|^2-4\langle a,v^0 \rangle.$

 \section{A general nonlinear selfdual variational principle for $\Lambda$-coercive functionals}
 
 In this section, we show that the ideas behind the nonlinear selfdual variational principles  can be extended in two different ways.  For one, and has already been noted in \cite{G3}, the hypothesis of regularity on the operator $\Lambda$ in Theorem \ref{main.1} can be weakened (see Definition \ref{pseudo} below).  We shall also relax the coercivity condition (\ref{one}) that  proved prohibitive in the case of evolution equations. 

  We start with the following weaker notion for regularity.
   
\begin{definition} \label{pseudo} \rm A map $\Lambda:D(\Lambda)\subset X\to X^*$ is said to be {\it pseudo-regular}  if  whenever $(x_n)_n$ is a sequence in $X$ such that $x_n \rightharpoonup x$ weakly in $X$ and $\limsup_n \langle \Lambda x_n, x_n-x\rangle \leq 0$, then 
$\liminf_n\langle \Lambda x_n, x_n\rangle\geq \langle \Lambda x, x\rangle$  and $\Lambda x_n \rightharpoonup \Lambda x$ weakly in $X^*$.
 \end{definition}
It is clear that regular operators are necessarily pseudo-regular operators. \\

 We also introduce the following weakened notion of coercivity. 

  \begin{definition}\rm  Let $J$ be the duality map from a reflexive Banach space $X$ into its dual  $X^*$, and consider a map $\Lambda:D(\Lambda)\subset X\to X^{*}$. Say that  a Lagrangian $L$  on $X \times X^*$   is {\it $\Lambda$-coercive}  if for any sequence $\{x_n\}_{n=1}^{\infty} \subseteq X$ such that $\|x_n\|\to +\infty$, we have 
\begin{equation}
\lim\limits_{\|x_n\|\to +\infty}L(x_n,\Lambda x_n+\frac{1}{n}Jx_n)+\langle x_n, \Lambda x_n \rangle +\frac {1}{n} \|x_n\|^2 =+\infty.
\end{equation}
 \end{definition}
The following is an extension of Theorem \ref{main.1}.
\begin{theorem} \label{main.10} Let $L$ be an  anti-selfdual Lagrangian on a reflexive Banach
space $X$ such that $0 \in  {\rm Dom}( L)$ and  ${\rm Dom}_1( L)$ is closed.   Let  $\Lambda: D(\Lambda)\subset X\to X^{*}$  be  a bounded pseudo-regular map 
such that ${\rm Dom}_1( L)\subset D(\Lambda)$, 
\begin{equation}
  \label{coercive}
\hbox{$L$ is $\Lambda$-coercive and $\langle\bar \partial L(x)+ \Lambda  x,  x\rangle \geq -C(\|x\|+1)$ for  large   $\|x\|$.}
\end{equation}
Then  there exists $\bar u \in X$ such that:
  \begin{eqnarray}  \label{principle}
 L(\bar u, \Lambda \bar u)+\langle \Lambda \bar u, \bar u\rangle &=&\inf_{x\in X} L(x, \Lambda x) +\langle \Lambda x, x\rangle= 0,
  \end{eqnarray}
and $\bar u$
 is a solution of the differential inclusion:  
 \begin{eqnarray}
  \hfill  -\Lambda \bar u &\in & \bar \partial L (\bar u). \label{eqn:zero1}
  \end{eqnarray}
\end{theorem}

\begin{remark} Theorem \ref{main.10} is an extension of Theorem \ref{main.1} which claims that (\ref{principle}) holds under the following coercivity assumption on $L$ and $\Lambda$. 
\begin{equation}\label{coercivity.0}
   \lim\limits_{\|x\|\to +\infty}H_{L}(0,-x)+\langle \Lambda x, x\rangle=+\infty.
 \end{equation}
 Indeed, in order to show that condition (\ref{coercivity.0}) is stronger than (\ref{coercive}), note that  for each $(x,p) \in X \times X^*$,
\begin{eqnarray*}
L(x,p)=\sup\{\langle y,p\rangle -H_L(x,y); y\in X\} \geq  -H_L(x,0) \geq & H_L(0,-x), 
\end{eqnarray*}
in such a way that if  
$\|x_n\|\rightarrow +\infty$, then 
\begin{equation*}
\lim\limits_{n \rightarrow +\infty}L(x_n,\Lambda x_n+\frac{1}{n}Jx_n)+\langle x_n, \Lambda x_n \rangle
+\frac {1}{n} \|x_n\|^2 \geq \lim\limits_{n \rightarrow +\infty} H_{L}(0,-x_n)+\langle \Lambda x_n, x_n\rangle =+\infty, 
\end{equation*}
from which follows that  $L$
 is $\Lambda$-coercive. Moreover, we have for  large   $\|x\|$, 
 \[
 \langle\bar \partial L(x)+\Lambda  x,  x\rangle = L(x, \bar \partial L(x))+\langle \Lambda x,x\rangle \geq H_L(0, -x)+\langle \Lambda x,x\rangle\geq  -C(\|x\|+1). 
\]
\end{remark}
For the proof of Theorem \ref{main.10}, we shall need the following lemma

\begin{lemma} \label{type M} Let $L$ be an anti-selfdual Lagrangian on a reflexive Banach space $X,$  let $\Lambda: D(\Lambda) \subseteq X\to X^{*}$ be  a pseudo-regular map 
 and let $F: D(F)\subseteq X\to X^{*}$  be  a regular map.  Assume $(x_n)_n$ is a sequence in $D(\Lambda)\cap D(F)$ such that $x_n \rightharpoonup x $  and $\Lambda x_n \rightharpoonup y$
 for some $x \in X$ and $y \in X^*$.  If
 $L(x_n, \Lambda x_n+F x_n)+\langle \Lambda x_n+F x_n, x_n \rangle =0$ for each $n\in \N$, 
then necessarily 
$ L(x, \Lambda x+F x)+\langle \Lambda x+F x, x \rangle 
 =0.$
\end{lemma}
\noindent{\bf Proof:} We have
\begin{eqnarray}\label{type M-1}
\limsup_{n}\langle \Lambda x_n, x_n-x \rangle &\leq & \lim_{n \rightarrow \infty}\langle \Lambda x_n,-x \rangle+\limsup_{n}
\big \{  -L(x_n, \Lambda x_n+F x_n)-\langle F x_n, x_n \rangle \big \}\nonumber \\
&=& \langle y,-x \rangle-\liminf_n \big \{  L(x_n, \Lambda x_n+F x_n)+\langle F x_n, x_n \rangle \big \}
\end{eqnarray}
Since $L$ is weakly lower semi continuous and $F$ is regular,  we have
\begin{eqnarray*}
L(x, y+F x)+\langle F x, x \rangle \leq \liminf_n  \big \{  L(x_n, \Lambda x_n+F x_n)+\langle F x_n, x_n \rangle \big \}
\end{eqnarray*}
which together with (\ref{type M-1}) imply
\begin{eqnarray*}
\limsup_{n}\langle \Lambda x_n, x_n-x \rangle & \leq&  \langle y,-x\rangle -L(x, y+F x)-\langle F x, x \rangle\\
&=& \langle y+F x,-x\rangle -L(x, y+F x).
\end{eqnarray*}
$L$ being an anti-selfdual Lagrangian, we have  $L(x, y+F x)\geq \langle y+ F x,- x \rangle,$ and therefore
\begin{eqnarray*}
\limsup_{n}\langle \Lambda x_n, x_n-x \rangle \leq  0.
\end{eqnarray*}
Now since $\Lambda$ is pseudo-regular,  
we have $y=\Lambda x$ and  $\liminf_{n}\langle \Lambda x_n, x_n \rangle \geq \langle \Lambda x, x \rangle. $
It follows from these facts that
\begin{eqnarray*}
 L(x, \Lambda x+F x)+\langle \Lambda x+F x, x \rangle \leq \liminf_{n} L(x_n, \Lambda x_n+F x_n)+\langle \Lambda x_n+F x_n, x_n \rangle =0,
  \end{eqnarray*}
On the other hand,  since $L$  is an anti-selfdual Lagrangian, we have the reverse inequality  $ L(x, \Lambda x+F
x)+\langle \Lambda x+F x, x \rangle \geq 0$ which implies the latter
to be  actually zero.$\Box$\\
 
\noindent{\bf Proof of Theorem \ref{main.10}:} Let $w(r)=\sup \{\|\Lambda u\|_*+1 ; \|u\|\leq r\}$ and let $Fu:=w(\|u\|) Ju.$
Let $L_\lambda^2$ be the $\lambda-$regularization of $L$ respect to the second variable  i.e.
\begin{eqnarray*}
L^{2}_\la (x,p):=\inf\left\{ L(x,q)+\frac{\| p-q\|_*^2}{2\la}+\frac{\la}{2}\| x\|^2; \, q\in X^*\right\} .
\end{eqnarray*}
Since $0 \in  {\rm Dom}( L)$  the Lagrangian $L$ and consequently $L_\lambda^2$ and therefore $H_{L_\lambda^2}(0,.)$
 are bounded from below. Also we have   $$ \lim\limits_{\|x\|\to +\infty}H_{L_\lambda^2}(0,-x)+\langle \Lambda x+\epsilon Fx, x\rangle=+\infty,$$
since $\langle \Lambda x+\epsilon Fx, x\rangle \geq -w(\|x\|)\|x\|+ \epsilon w(\|x\|)\|x\|^2.$

It follows from  Theorem \ref{main.1} that there exists
$x_{\epsilon, \lambda}$ such that
\begin{eqnarray*}
 L_\lambda^2(x_{\epsilon, \lambda}, \Lambda x_{\epsilon, \lambda}+\epsilon Fx_{\epsilon, \lambda})+\langle \Lambda x_{\epsilon, \lambda}+\epsilon Fx_{\epsilon, \lambda}, x_{\epsilon, \lambda}\rangle =0
  \end{eqnarray*}
  which means that $\Lambda x_\epsilon +\epsilon Fx_\epsilon  \in -\bar\partial L_\lambda^2 (x_\epsilon)$, and in other words,
$\Lambda x_{\epsilon, \lambda}+\epsilon Fx_{\epsilon, \lambda}
+\lambda Jx_{\epsilon, \lambda} \in -\bar \partial L(x_{\epsilon, \lambda})$.  This
together with (\ref{coercive}), imply
$\langle \epsilon Fx_{\epsilon, \lambda} +\lambda Jx_{\epsilon, \lambda}, x_{\epsilon, \lambda}\rangle \leq C \|x_{\epsilon, \lambda}\|,$
thereby giving
\begin{eqnarray*}
 \epsilon w(\|x_{\epsilon, \lambda}\|) \|x_{\epsilon, \lambda}\|^2+\lambda  \|x_{\epsilon, \lambda}\|^2 \leq C \|x_{\epsilon, \lambda}\|,
\end{eqnarray*}
which in turn  implies that  $(\epsilon F x_{\epsilon, \lambda} )_\epsilon$ and $ (x_{\epsilon, \lambda})_\epsilon$ are bounded.  Since now $\Lambda$ is a bounded operator, we get that $\Lambda  x_{\epsilon, \lambda}$ is bounded in $X^*$.   Suppose, up to a subsequence, $x_{\epsilon, \lambda}\rightharpoonup x_{\lambda}$
 and $\Lambda  x_{\epsilon, \lambda}\rightharpoonup p_\lambda.$  It follows from Lemma \ref{type M} that for every $\lambda>0$, we have
\begin{eqnarray*}
 L(x_{\lambda}, \Lambda x_{\lambda}+\lambda Jx_{\lambda})+\langle \Lambda x_{\lambda}+\lambda Jx_{\lambda}, x_{\lambda}\rangle = 0.
  \end{eqnarray*}
Since $L$ is $\Lambda$-coercive, $x_{\lambda}$ is a bounded sequence in $X$ and therefore converges  weakly -- up to a subsequence-- to a $\bar u \in X.$
 Again, since $\Lambda$ is a bounded operator, $\Lambda x_{\lambda}$
 is also bounded in  $X^*,$  and again Lemma \ref{type M} yields
$L(\bar u, \Lambda \bar u)+\langle \Lambda \bar u, \bar u\rangle = 0$,
which means that $-\Lambda  \bar u \in \bar \partial L(\bar u ).$
\begin{remark} Note that, we do not really need  that $\Lambda$ is a bounded operator, but a weaker  condition of the form $\|\Lambda x\| \leq C H (0,x)+ w(\|x\|)$  for some nondecreasing function $w$ and some constant $C>0$. 
\end{remark}
\hfill $\square$\\

Let now $A: D(A)\subset X\to X^*$  be a closed linear 
operator on a reflexive Banach space $X$, and consider  $X_A$ to be the Banach space that is the closure of $D(A)$ 
for the norm $\|x\|_A=\|x\|_X+\|Ax\|_{X^*}$.  We have the following result. 
\begin{corollary}\label{cor.1} Let $A: D(A)\subset X\to X^*$  be a closed linear 
operator on a reflexive Banach space $X$ with a dense domain, and   let $\Lambda $ be a map from $D(A)$  into  $X^*$ that induces a pseudo-regular operator  $\Lambda: X_A\rightarrow X_A^*$.
Suppose $L$ is an anti-selfdual Lagrangian on $X \times X^*$ that satisfies the following conditions:
\begin{equation}
\hbox{ $L$ is $(\Lambda+A)$-coercive on $X_A$}  
\end{equation}
\begin{equation}
\hbox{For each $p \in {\rm Dom}_2(L)$, the functional  $x\rightarrow L(x,p)$ is continuous on $X$.}
\end{equation}
\begin{equation}
\hbox{ $x\to L(x,0)$ is bounded  on the unit ball of $X$.}
\end{equation}
Then  there exists $\bar u \in X_A$ such that:
  \begin{eqnarray*}
 L(\bar u, \Lambda \bar u+A \bar u)+\langle \Lambda \bar u+A \bar u, \bar u\rangle &=&\inf_{x\in Y} L(x, \Lambda x+A x) +\langle \Lambda x+A x, x\rangle= 0,\\
 -\Lambda \bar u-A \bar u &\in&  \bar \partial L (\bar u).
  \end{eqnarray*}
\end{corollary}

\noindent{\bf Proof:} Note first that  $X_A \subseteq X \subseteq X^* \subseteq X_A^*. $ We first show that the Lagrangian 
\begin{eqnarray*}
{\cal M}(u,p):= \left\{ \begin{array}{lcl}
 L(u,p), \quad  p \in  X^*\\
  +\infty  \qquad     p \in  X_A^*\setminus X^*
\end{array}\right.
\end{eqnarray*}
is an anti-selfdual Lagrangian on $X_A \times X_A^*$.  Indeed,  if $q \in X^*$, use the fact that $X_A$ is dense in $X$ and that  the functional  $x\rightarrow L(x,p)$ is continuous on $X$ to write
\begin{eqnarray*}
{\cal M}^*(q,v)&=&\sup \{ \langle u, q\rangle+ \langle v, p\rangle-{\cal M}(u,p); (u,p) \in X_A \times X_A^* \}\\
&=&\sup \{ \langle u, q\rangle+ \langle v, p\rangle- L(u,p); (u,p) \in X_A \times X^* \}\\ &=&  L^*(q,v)=L(-v,-q)={\cal M}(-v,-q).
\end{eqnarray*}
If now $q \in  X_A^*\setminus X^*$, then  there exists $\{x_n\}_n \subseteq X_A $ with $\|x_n\|_X \leq 1$ such that
$\langle x_n, q \rangle \rightarrow +\infty  \text{  as  } n\rightarrow \infty.$
It follows,
\begin{eqnarray*}
{\cal M}^*(q,v)&=&\sup \{ \langle u, q\rangle+ \langle v, p\rangle-{\cal M}(u,p); (u,p) \in X_A \times X^* \}\\
&\geq &\sup \{ \langle x_n, q\rangle- L(x_n,0) \}\\ &=&+\infty={\cal M}(-v,-q)
\end{eqnarray*}
Note that since  $\{x_n\}_n$ is a bounded sequence in $X,$ the
sequence $\{L(x_n,0)\}_n$ is bounded. It follows from the
assumptions that ${\cal M}$ is $\Lambda+A$ coercive on $X_A$  and the rest of condition (\ref{coercive}) in Theorem \ref{main.10}. $\square$\\
\begin{corollary} \label{cor.2} Let the operator $A$ and the space $X_A$ be as in Corollary \ref{cor.1} and  let $\phi$ be a proper convex lower semi-continuous function that is both coercive and  bounded in $X$. Let   $\Lambda:  X_A\to X^{*}$ be a pseudo-regular operator and assume the following conditions: 
 \begin{equation}\label{1}
\hbox{  $u \rightarrow \langle u, \Lambda u+ A u \rangle $ is  bounded from below.}
    \end{equation}
    \begin{equation}\label{k-condition}
 \hbox{   $\|\Lambda u\|_{X^*}\leq k \|A u\|_{X^*}+ w (\|u\|_X) $ for some constant $0< k<1$ and a nondecreasing function $w$.}
  \end{equation}
   Then  a solution $\bar x \in X_A$ to the equation 
$ 0\in  \Lambda x+Ax+\partial \varphi (x)$
 can be obtained as a minimizer of the problem:
\begin{equation*}
 \label{min.99}
\inf_{x\in X}\left\{\phi (x)  +\phi^{*}(-\Lambda x- Ax) +\langle x, \Lambda x +Ax\rangle \right\}=0.  
\end{equation*}
   \end{corollary}
\noindent {\bf Proof:} It is an immediate consequence of Corollary \ref{cor.1} applied to the Lagrangian $L(x,p)=\phi(x)+\phi^*(-p).$ To verify the conditions of \ref{cor.1}, we note that 
\begin{eqnarray*}
\langle  \partial \phi (u)+ \Lambda u+ A u, u \rangle \geq \langle  \partial \phi (0), u \rangle+\langle   \Lambda u+ A u, u \rangle \geq -C(1+\|u\|_X) \geq -C(1+\|u\|_Y).
\end{eqnarray*}
We have used the fact that $\phi $ is convex monotone and $\langle   \Lambda u+ A u, u \rangle $ is bounded from below. Now we prove that $L$ is $(\Lambda +A)$ coercive. Indeed,  suppose $\{x_n\}_n\subseteq X_A$ is a subsequence such that $\|x_n\|_{X_A}\rightarrow \infty,$ we show that 
\begin{eqnarray*}
 \left\{\phi (x_n)  +\phi^{*}(-\Lambda x_n- Ax_n- \frac {1}{n}J x_n) +\langle x_n, \Lambda x_n +Ax_n\rangle +\frac {1}{n} 
\|x_n\|_X\right\} \rightarrow \infty.
\end{eqnarray*}
Indeed, if the above relation does not hold since $\langle x_n, \Lambda x_n +Ax_n\rangle +\frac {1}{n} 
\|x_n\|_X$ is bounded from below, we have $\phi (x_n)  +\phi^{*}(-\Lambda x_n- Ax_n- \frac {1}{n}J x_n) $ is bounded from above. The coerciveness of  $\phi$  on $X$ ensures the boundedness of $\{\|x_n\|_X\}_n$.  Now we show that $\{x_n\}$ is actually bounded in $X_A$. In fact , since $\phi$ is bounded  on $X$ we have that $\phi^*$ is coercive in $X^*$ and in result
\begin{eqnarray*}
\|\Lambda x_n+Ax_n+ \frac {1}{n}J x_n\|_{X^*}\leq  C
\end{eqnarray*}
for some constant $C>0.$  It follows from (\ref{1})   and the above that  
\begin{eqnarray*}
\|A x_n\|_{X^*} &\leq & \|\Lambda x_n+Ax_n+ \frac {1}{n}J x_n\|_{X^*}+\|\Lambda x_n+ \frac {1}{n}J x_n\|_{X^*}\\
&\leq& C+ \|\Lambda x_n\|_{X^*}+ \frac {1}{n}\|J x_n\|_{X^*}\\ &\leq & C+ k \|A x_n\|_{X^*} + w(\|x_n\|_X)+ \frac {1}{n}\| x_n\|_X
\end{eqnarray*}
Hence 
$(1-k)\|A x_n\|_{X^*} \leq C + w(\|x_n\|_X)+ \frac {1}{n}\| x_n\|_X$, 
and therefore $\|A x_n\|_{X^*}$ is bounded which results the boundedness of $\{x_n\}$ in $X_A.$ \hfill $\square$\\

We can also give a variational resolution for certain nonlinear  systems. 

 \begin{corollary} \label{cor.system} Let $\phi$ be a  bounded  convex lower semi-continuous function  on $X_1\times X_2$, let $A:X_1\to X_2^{*}$ be any bounded linear operator, let $B_{1}:X_1 \to X_1^{*}$ (resp., $B_{2}:X_2\to X_2^{*}$) be two positive linear operators. Let $Y_i:= \{x \in X_i ; B_i x \in X_i^*\}, i=1,2.$  Assume $\Lambda:=(\Lambda_{1}, \Lambda_{2}): Y_1\times Y_2 \to Y_1^{*}\times Y_2^{*}$ is a pseudo-regular operator such that  
 \[
 \lim\limits_{\|x\|_{X_1}+\|y\|_{X_2}\to \infty}\frac{\phi(x,y)+ \langle B_{1}x, x\rangle+ \langle B_{2}y, y\rangle+ \langle \Lambda (x,y), (x,y)\rangle }{\|x\|_{X_1}+\|y\|_{X_2}}=+\infty,
\]
and 
\[\|(\Lambda_1,\Lambda_2) (x,y)\|_{X_1^* \times X^*_2}\leq k \|(B_1,B_2) (x,y)\|_{X_1^* \times X^*_2}+ w(\| (x,y)\|_{X_1 \times X_2})\]
for some continuous and non-decreasing function $w$, and some constant $0<k<1$.   Then for any $(f, g)\in Y_1^{*}\times Y_2^{*}$, there exists $(\bar x,\bar y) \in Y_1\times Y_2$  which solves the following system
     \begin{equation*}
 \left\{ \begin{array}{lcl}
\label{eqn:existence}
\hfill -\Lambda_{1}(x,y)-A^{*} y-B_{1} x +f&\in& \partial_{1} \phi ( x,  y).\\
\hfill -\Lambda_{2}(x,y)+A x-B_{2} y+g&\in& \partial_{2} \phi ( x,  y).
\end{array}\right.
 \end{equation*}
 
The solution is obtained as a minimizer on $Y_1\times Y_2$ of the functional 
 \[
I(x,y)=\psi (x, y)+\psi^{*}(-A^{*}y-B_{1}x-\Lambda_{1}(x,y), Ax-B_{2}y-\Lambda_{2}(x,y))+ \langle B_{1}x, x\rangle+ \langle B_{2}y, y\rangle+ \langle \Lambda (x,y), (x,y)\rangle.
\]
where 
\[
\psi (x,y)=\phi (x,y) -\langle f, x\rangle- \langle g, y\rangle.
\]

 \end{corollary} 
\noindent {\bf Proof:}  Consider the following ASD Lagrangian  (see \cite{G2})
  \[
L((x,y), (p,q))=\psi(x, y)+\psi^{*}(-A^{*}y-p, Ax-q).
\]
Setting $B:= (B_1, B_2)$,  Corollary  \ref{cor.2} yields that $I(x,y)=L((x,y), \Lambda (x,y)+B(x,y))+ \langle\Lambda (x,y)+B(x,y) , (x,y)\rangle $ attains its minimum at some point $(\bar x,\bar y) \in Y_1\times Y_2$ and that the minimum is $0$. In other words, 
\begin{eqnarray*}
0&=&I(\bar x,\bar y)\\&=&
\psi(\bar x, \bar y)+\psi^{*}(-A^{*}\bar y-B_{1}\bar x-\Lambda_{1}(\bar x,\bar y), A\bar x-B_{2}\bar y-\Lambda_{2}(\bar x,\bar y))+ \langle\Lambda (\bar x, \bar y)+B(\bar x, \bar y) , (\bar x, \bar y)\rangle  \\
&=&\psi (\bar x, \bar y)+\psi^{*}(-A^{*}\bar y-B_{1}\bar x-\Lambda_{1}(\bar x,\bar y), A\bar x-B_{2}\bar y-\Lambda_{2}(\bar x,\bar y))+\\
&&\langle(\Lambda_1 (\bar x, \bar y)+B_1\bar x -A^* \bar y, \Lambda_2 (\bar x, \bar y)+ B_2 \bar y +A \bar x) , (\bar x, \bar y)\rangle
\end{eqnarray*}
from which  follows that
\begin{equation*}
 \left\{ \begin{array}{lcl}
\label{eqn:existence}
-A^{*}y-B_{1}x-\Lambda_{1}(x,y)  &\in& \partial_{1} \varphi (x,y)-f\\
\hfill Ax-B_{2}y-\Lambda_{2}(x,y)&\in& \partial_{2} \varphi (x, y) -g.
\end{array}\right.
\end{equation*}

 \subsection*{Example : A variational resolution for doubly nonlinear coupled equations}
  
Let ${\bf b_{1}}:\Omega \to  {\bf R^{n}}$ and  ${\bf b_{2}}:\Omega \to  {\bf R^{n}}$ be two compact supported smooth vector fields on the neighborhood of a bounded domain $\Omega$ of  $\bf R^{n}$.    Consider the Dirichlet problem:
\begin{equation}
\label{Ex1.500}
 \left\{ \begin{array}{lcl}
    \hfill Ê\Delta v +{\bf b_{1}}\cdot \nabla u &=& |u|^{p-2}u+u^{m-1}v^{m} +f  \hbox{\rm \, on \, $\Omega$}
\\
  \hfill -Ê\Delta u +{\bf b_{2}}\cdot \nabla v&=& |v|^{p-2}v -u^{m}v^{m-1} +g  \, \hbox{\rm \, on \, $\Omega$}\\
 \hfill  u=v &=& 0 \quad \quad \quad \quad \quad \quad \quad  \quad \quad \quad \hbox{\rm on  $\partial \Omega$. }        \end{array}  \right.
   \end{equation}
 We can use Corollary \ref{cor.2} to get
 \begin{theorem} Assume  $f,g$ in $L^p$, $2 \leq p$, that   ${\rm div }({\bf b_{1}})\geq 0$ and ${\rm div }({\bf b_{2}})\geq 0$ on $\Omega$, and $1\leq m<\frac{p-1}{2}$. Let $X= \{u \in H^{1}_{0}(\Omega); u \in L^p(\Omega)  \&  Ê\Delta u \in L^q(\Omega)\}$  and consider on $X \times X $ 
 the functional 
\begin{eqnarray*}
I(u, v)&=&\Psi (u) +\Psi^{*}({\bf b_{1}}.\nabla u +\Delta v-u^{m-1}v^{m} )+\Phi (v) +\Phi^{*}({\bf b_{2}}.\nabla v  -\Delta u +u^{m}v^{m-1})\\
&&+\frac{1}{2}\int_{\Omega}{\rm div }({\bf b_{1}}) \, |u|^{2} dx+\frac{1}{2}\int_{\Omega}{\rm div }({\bf b_{1}}) \, |v|^{2} dx
\end{eqnarray*}
where
\[
\hbox{$ \Psi (u)=\frac{1}{p}\int_{\Omega}|u|^{p} dx+\int_{\Omega}fu dx$, and $
 \Phi (v)=\frac{1}{p}\int_{\Omega}|v|^{p} dx+\int_{\Omega}gv dx$} 
 \]
are defined on $L^p(\Omega)$ and $\Psi^{*}$ and $\Phi^{*}$ are their Legendre transforms in $L^q(\Omega)$. Then there exists $(\bar u, \bar v)\in X \times X$ such that:
\[
I(\bar u, \bar v)=\inf \{I(u, v);  (u, v)\in X\times X\}=0, 
\]
and $(\bar u, \bar v)$ is a solution of $(\ref{Ex1.500})$.
 \end{theorem}
\noindent{\bf Proof:} Let $A=\Delta$, $X_A= X$ and $X_1=L^p(\Omega).$  $\Phi$ and $\Psi$ are continuous and coercive on $X_1.$  To show that $I$ is $\Lambda-$coercive, we just need to verify condition (\ref{k-condition}) in Corollary \ref{cor.2}. Indeed, by H\"older's inequality for $q= \frac{p}{p-1}\leq 2$ we obtain
\begin{eqnarray*}
\|u^{m}v^{m-1}\|_{L^q(\Omega)}  \leq  \|u\|^{m}_{L^{2mq}(\Omega)}\|v\|^{(m-1)}_{L^{2(m-1)q}(\Omega)}
 \end{eqnarray*}
and since $m < \frac {p-1}{2}$ we have $2mq < p$ and therefore 
\begin{eqnarray}\label{k-condition1}
\|u^{m}v^{m-1}\|_{L^q(\Omega)}  \leq  C( \|u\|^{2m}_{L^{p}(\Omega)}+\|v\|^{2(m-1)}_{L^{p}(\Omega)}.
 \end{eqnarray}
Also  since $q\leq 2$,
\begin{eqnarray}\label{k-condition2}
\|{\bf b_{1}}.\nabla u\|_{L^q(\Omega)} & \leq & C\|{\bf b_{1}}\|_{L^{\infty}(\Omega)} \|\nabla u\|_{L^2(\Omega)} \nonumber \\& \leq & C\|{\bf b_{1}}\|_{L^{\infty}(\Omega)}\big (\int \langle -\Delta u, u \rangle  \, dx \big )^{\frac{1}{2}}\leq C\|{\bf b_{1}}\|_{L^{\infty}(\Omega)} \|u\|^{\frac{1}{2}}_{L^p(\Omega)}\|\Delta u\|^{\frac{1}{2}}_{L^q(\Omega)} \nonumber \\& \leq & k\|\Delta u\|_{L^q(\Omega)}+ C(k)\|{\bf b_{1}}\|^2_{L^{\infty}(\Omega)} \|u\|_{L^p(\Omega)}
 \end{eqnarray}
for some $0<k<1.$  Hence  condition (\ref{k-condition}) follows from (\ref{k-condition1}) and (\ref{k-condition2}).
 
Also, it is also easy to verify that the nonlinear operator $\Lambda: X \times X \to L^q(\Omega) \times L^q(\Omega)$ defined by
\[
\Lambda (u,v)=(-u^{m-1}v^{m}+{\bf b_{1}}.\nabla u, u^{m}v^{m-1}+{\bf b_{2}}.\nabla v)
\]
is regular. \hfill $\square$\\

It is worth noting that there is no restriction on the power $p$ in the previous example, that is $p$ can well be beyond the critical Sobolev  exponent.

\end{document}